\documentclass[12pt,draft]{article}

\usepackage{amsfonts}
\usepackage{graphics}
\usepackage{amsmath}
\usepackage{amscd}
\usepackage{amssymb}
\usepackage{euscript}

\newtheorem{Theorem}{Theorem}[section]

\newtheorem{Corollary}[Theorem]{Corollary}
\newtheorem{Example}[Theorem]{Example}
\newtheorem{Definition}[Theorem]{Definition}

\newtheorem{Remark}[Theorem]{Remark}

\newtheorem{Lemma}[Theorem]{Lemma}
\newtheorem{Proposition}[Theorem]{Proposition}

\newtheorem{Fundamental Theorem}{Fundamental Theorem}

\newenvironment{Proof}[1][Proof]{\textbf{#1.} }{\ \rule{0.5em}{0.5em}}

\def\leq{\leqslant} 
 
\def \ge {\geqslant}

\def \A {\mathcal{A}}
\def \B {\mathcal{B}}

\def \C {\mathcal{C}}
\def \CT {\mathrm{cotrunc}}
\def \coker {\mathrm{coker}}
\def \d {\partial}
\def \D {{ \cal D}}

\def \f {\phi}

\def \Fc {\mathcal{F}}

\def \g {\gamma}
\def \G {\mathcal{G}}
\def \Hom {\mathrm{Hom}}

\def \id {\mathrm{id}}
\def \il {{^L_{i_L}}}

\def \l {\lambda}
\def \L {\Lambda}
\def \M {\hat{M}}
\def \Mc {\mathcal{M}}
\def \Mce {\EuScript{M}}
\def \N {\mathbb{N}}
\def \Nc {\mathcal{N}}
\def \p {\psi}
\def \r {\rho}
\def \ra {\xrightarrow}

\def \S  {\Sigma}
\def \t {\triangleright}

\def \V {\bigvee}

\def \w {\omega}
\def \X {\EuScript{X}^\pi}
\def \Z {\mathbb{Z}}

\begin{document}

\title{On the Homotopy Type and the Fundamental Crossed Complex of the Skeletal Filtration of a CW-Complex}

\author{ Jo\~{a}o  Faria Martins\\ \footnotesize\it  {Departamento de Matem\'{a}tica, Instituto Superior T\'{e}cnico,}\\ {\footnotesize\it Av. Rovisco Pais, 1049-001 Lisboa, Portugal}\\ {\footnotesize\it jmartins@math.ist.utl.pt}}

\date{\today}

\maketitle

\begin{abstract}
We prove that if $M$ is a CW-complex, then the homotopy type of the skeletal
filtration of $M$ does not depend on the cell decomposition of $M$ up to wedge
products with $n$-disks $D^n$, when the later are given their natural
CW-decomposition with unique cells of order $0$, $(n-1)$ and $n$; a result
resembling J.H.C. Whitehead's work on simple homotopy types. From the Colimit
Theorem for the  Fundamental Crossed Complex of a CW-complex (due to R. Brown and P.J. Higgins),  follows an algebraic analogue for the fundamental
crossed complex $\Pi(M)$ of the skeletal filtration of $M$, which thus 
depends only on the homotopy type of $M$ (as a space) up to free product with
crossed complexes of the type $\D^n \doteq \Pi(D^n), n \in \N$. This expands an old result
(due to J.H.C. Whitehead) asserting that the homotopy type of $\Pi(M)$
depends only on the homotopy type of $M$.
We use these results to define a homotopy invariant $I_\A$ of CW-complexes for any
finite crossed complex $\A$. We  interpret it in terms of the weak homotopy
type of the  function space $TOP((M,*),(|\A|,*))$, where $|\A|$ is the
classifying space of the crossed complex $\A$. 
\end{abstract}

2000 Mathematics Subject Classification: 55P10, 55Q05, 57M05.

\tableofcontents

\section*{Introduction}

The concept of a ``crossed complex'', a type of group complex, was invented by
A. L. Blackers under the name of ``Group System''. See \cite{Bl}.  In
\cite{W4,W5}, J.H.C Whitehead, re-introduced them,
considering, however, only the important totally free
case.  Crossed complexes,  algebraically describe, for instance, the complex of groups:
\[... \ra{\d_5} \pi_4 (M^4,M^3,*) \ra{\d_4} \pi_3 (M^3,M^2,*) \ra{\d_3} \pi_2(M^2,M^1,*) \ra{\d_2}
\pi_1(M^1,*),\]
with  the obvious boundary maps $\d_n, n \in \N$, together  with the standard action of $\pi_1(M^1,*)$ on $\pi_n(M^n,M^{n-1},*),n>1$. Here $M$ is a connected CW-complex and $M^n$ is its  $n$-skeleton.   This crossed complex, denoted by $\Pi(M)$,  is
called the Fundamental Crossed Complex of $M$ (provided with its skeletal
filtration). More generally, we can take $M$ to be any filtered space.  What distinguishes crossed complexes from the usual group complexes, is that they  have an action of the first group on all the others, satisfying some natural compatibility conditions, which are stronger than merely imposing that the boundary maps preserve the actions.

Crossed complexes were extensively studied and used, for example, by
J.H.C. Whitehead, H.J.  Baues, R. Brown and P.J. Higgins, see
\cite{W4,W5,BA1,BA2,B1,BH2,BH3,BH4,BH5,BS}. They admit an obvious notion of
homotopy, as well as classifying spaces. Another useful feature is that the
category of crossed complexes is a category with colimits, thus in particular
the free product of any family of crossed complexes is well defined.  In
addition,  there
exists a general  ``van Kampen type property'', so that the fundamental
crossed complex functor from the category of CW-complexes to the category of
crossed complexes preserves colimits. This result is due to R. Brown and P.J. Higgins, and appeared in \cite{BH3,BH4} under the designation  ``General van Kampen Theorem''.

Let $M$ be any space which can be given the structure of a CW-complex. The
fundamental crossed complex $\Pi(M)$ of the skeletal filtration of $M$ is strongly dependent on the chosen
cellular structure. However, J.H.C.  
Whitehead proved (see \cite{W4}) that if $M$ is a CW-complex, then the
homotopy type of the crossed complex $\Pi(M)$ depends only on the homotopy
type of $M$ as a space.  In fact, the results of J.H.C. Whitehead that  appeared in
\cite{W4,W5} imply immediately  a stronger result: If
$M$ and $N$ are homotopic CW-complexes of dimension $\leq n$, then $\Pi(M')$
and $\Pi(N')$  are simply homotopy equivalent. Here $M'$ and $N'$ are
obtained from $M$ and $N$ by taking wedge products with a certain number of
spheres $S^n$. We refer to \cite{W5} for the definition of simple homotopy
equivalence of crossed complexes.

In this article we prove a  theorem expanding the first of the results due to
J.H.C. Whitehead   which we referred to in the previous paragraph; a result  in the direction of the second of them.  Namely, we show that  $\Pi(M)$
depends only on the  homotopy type of $M$, as a space, up to free products
with crossed complexes of the type $\D^n \doteq \Pi(D^n)$, where $n \in\N$. Here $D^n$ is
the $n$-disk with its natural CW-decomposition with one 0-cell, one
$(n-1)$-cell and one $n$-cell. This theorem is deduced (making use of the
General van Kampen Theorem) from  an analogous statement on CW-complexes (proved
in this article): The homotopy type of a CW-complex $M$, as a filtered space,
depends only on the homotopy type of $M$ (as a space), up to wedge products
with CW-complexes of the type $D^n$, where $n\in\N$, provided with their natural
cell decompositions, already described. This has strong similarities with the following result of
J.H.C. Whitehead on 
simple homotopy types of CW-complexes: If
$M$ and $N$ are CW-complexes of dimension $\leq n$, and $f:M \to N$ is a
$(n-1)$-equivalence, then $M'$ and $N'$  are simply homotopy equivalent, where $M'$ and $N'$ are
obtained from $M$ and $N$ by taking wedge products with a certain number of
spheres $S^n$. Note that  the previously stated result on the simple homotopy
types of $\Pi(M)$ and $\Pi(N)$ can be  deduced from this statement  together
with  the General van Kampen Theorem.

As an application, we prove that if $\A$ is a finite crossed complex, then the
number of morphisms from $\Pi(M)$ into $\A$ can be naturally normalised to a
homotopy invariant $I_\A(M)$ of CW-complexes $M$, which is easy to calculate,
in fact as simply as the cellular homology groups of $M$. This invariant generalises
a previous construction for crossed modules appearing in \cite{FM2}. There
it was proved  that, if $\A$ is a crossed module, then  the value of $I_\A$
on  a knot complement defines a non-trivial  invariant of knotted
surfaces $\S$. Moreover, we elucidated a graphical  algorithm for its calculation
from a movie representation of $\S$.

The homotopy invariant $I_\A$  also upgrades previous constructions of invariants of
manifolds and knots derived from finite crossed modules, see
\cite{Y,P1,FM1}. Another construction by T. Porter appearing in \cite{P2} (from which T.Q.F.T.'s can be defined) has
the merit of considering  finite $n-cat$ groups, which are more general than
finite 
crossed complexes. Finally, in \cite{Mk},  M. Mackaay
defined  4-manifold invariants, conjecturally related to the $n=3$ case of
T. Porter's construction, but considering, furthermore, an additional twisting
by cohomology classes  of $3$-types. This twisting is  similar to the one R. Dijkgraaf and
E. Witten's introduced in \cite{DW}, even though their 3-dimensional oriented manifold
invariant (the well known 
Dijkgraaf-Witten invariant)  
considers only finite $1$-types (in other words, finite groups).

At the end of this article, we interpret and give an alternative proof of the
existence of the invariant $I_\A$, where $\A$ is a finite crossed
complex. Explicitly,  we prove that $I_\A(M)$ is determined by the number of
connected components as well as   the homotopy groups of the function space
$TOP((M,*),(|\A|,*))$, where $|\A|$ is the classifying space of $\A$; by a
multiplicative Euler Characteristic type formula. This
result is  an application of the general theory of classifying spaces of  crossed
complexes developed by R. Brown and P.J. Higgins, see \cite{BH5}, and also of the non-trivial fact that given two crossed complexes $\A$ and $\B$ we can define a crossed complex $CRS(\A,\B)$, made from all of crossed complex morphisms  $\A \to \B$, together with  their $n$-fold homotopies, where $n \in \N$. See \cite{BH3,BH4,BH5}.

This description of $I_\A$ enables us to incorporate $n$-dimensional
cohomology classes of $|\A|$ into it, as long as we  consider only
$n$-dimensional closed oriented manifolds. This yields therefore an extension of
Dijkgraaf-Witten's invariant to crossed complexes. Its full description
will appear in a future paper.

\section{Preliminaries}
All CW-complexes considered in this article will  be connected, with an extra
 technical condition imposing that they have a unique $0$-cell, taken as their
 base point. We will very often make the assumption that they only have a
 finite number of $n$-cells for any $n \in \N$. If $M$ is a CW-complex, denote
 the $n$-skeleton of $M$ by $M^n$, where $n \in \N$. If $f:M \to N$ is a
 cellular map, where $N$ is a CW-complex, denote $f^n\doteq f_{|M^n}: M^n \to
 N^n$. An $n$-type is a CW-complex $M$ such that $\pi_k(M)=\{0\}$ if $k>n$.
 Set $I=[0,1]$.
\subsection{Cofibred Filtrations}
Most issues treated in this subsection are widely known. The specialist should
pass directly to   section 2, and use this subsection, as well as the next one,  only as a reference. 

\begin{Definition}
Let $M$ be  some (Hausdorff) space. A filtration of $M$ is a sequence $\{M_n\}_{n =0}^\infty$ (briefly $\{M_n\}$) of subspaces of $M$ such that $M_n \subset M_{n+1}, \forall n \in \N$, and also $M= \bigcup_{n \in \N} M_n$. A space $M$ provided with a filtration is called a filtered  space.  A filtration  of $M$ is called finite, of length $L\in \N$, if $M_n=M, \forall n\ge L$.  A filtration $\{M_n\}_{n \in \N}$ of $M$ is called cofibred if the inclusion   $M_n \to M_{n+1}$ is a cofibration for any $n \in \N$. If $M$ and $N$ are filtered topological spaces, then a filtered map $f\colon M\to N$ is a continuous map such that $f(M_n) \subset N_n, \forall n \in \N$.
\end{Definition}
Obviously filtered topological spaces and filtered maps form a category. If
$f\colon M \to N$ is a filtered map, we set $f_n\doteq f_{|M_n}\colon M_n \to N_n$, where
$ n \in \N$. Considering the relation of filtered homotopy, we can consider
the category of filtered topological spaces, with morphisms the filtered maps
up to  filtered homotopy equivalence.

\begin{Example}
Let $M$ be a CW-complex. Then $M$ has a natural cofibred filtration
$\{M^n\}_{n \in \N}$, where $M^n$ is the n-skeleton of $M$, the ``Skeletal
Filtration of $M$''. More generally,  if $\{M_n\}_{n \in \N}$ is a filtration of $M$, such that each $M_n$ is a CW-complex included in $M_{n+1}$, cellularly, for any $n \in \N$, then   $\{M_n\}_{n \in \N}$ is a cofibred filtration. In fact, all the filtrations that we consider in this article will be of this particular type, called ``Filtrations by Subcomplexes''.
\end{Example}

\subsubsection{Some   Results on Cofibre Homotopy Equivalence}\label{important}
The following result is well known:
\begin{Theorem}\label{Simplify}
Let $M$ and $N$ be spaces provided with cofibred filtrations, say  $\{M_n\}_{n \in \N}$ and $\{N_n\}_{n \in \N}$, which we suppose to be finite (we will eliminate this condition later, in the cellular case). Let $F\colon M \to N$ be a filtered map. Suppose that each map $F_n\colon M_n \to N_n, n \in \N$ is a  homotopy equivalence. Then it follows that $F$ is a filtered homotopy equivalence. 
\end{Theorem}
 This is shown in \cite{M} (end of chapter 6), or \cite{B2} (7.4), for the case in which both filtrations of $M$ and $N$ have length $L=1$. The proof done in \cite{B2} extends immediately to the general case, inductively, as we will  indicate in a second.  We will give later  an alternative proof of Theorem \ref{Simplify} for the particular case of filtrations by subcomplexes, which adapts to infinite filtrations, a generality we will need to be able to work with infinite dimensional CW-complexes.

It is convenient to recall the following result appearing in \cite{B2} (7.4), which leads immediately  to a proof of the  Theorem \ref{Simplify}:

\begin{Lemma} {\bf (R. Brown)}
Let $F\colon(M_1,M_0) \to (N_1,N_0)$ be a map between cofibred pairs. Suppose that $F=F_1\colon M_1 \to N_1$ and $F_0\colon M_0 \to N_0$ are homotopy equivalences. Let $G_0\colon N_0 \to M_0$ be a homotopy inverse of $F_0$, and consider homotopies (in $M_0$ and $N_0$):
\[K\colon G_0F_0 \to \id_{M_0} \textrm{ and } H\colon F_0G_0 \to \id_{N_0}.\]
Then $G_0$ extends to a homotopy inverse $G$ of $F$, such that there exists not only a  homotopy $FG \to \id_{N_1}$ which extends $H$, but also  a  homotopy  $GF \to \id_{M_1}$ extending the concatenation $J$ of the following homotopies:
\[G_0F_0 \ra{G_0F_0 K^{-1}} G_0 F_0 G_0  F_0\ra {G_0 H F_0} G_0 F_0 \ra{K} \id_{M_0}.\]
 Note that if $G_0,F_0,H,K$ are filtered, for some filtrations of  $M_0$ and $N_0$, then so is $J$.
\end{Lemma}
Therefore, Theorem \ref{Simplify} can be strengthened in the following way:

\begin{Corollary}\label{strenght}
Under the conditions of Theorem \ref{Simplify}. If: 
\[G_L\colon(N_L,N_{L-1},...,N_0) \to (M_L,M_{L-1},...,M_0)\]
is a filtered homotopy inverse  of the filtered map:  
\[F_L=F_{|M_L} \colon (M_L,M_{L-1},...,M_0) \to (N_L,N_{L-1},...,N_0),\]
 then there exists a filtered homotopy inverse  $G$ of $F$ extending $G_L$.
\end{Corollary}

Now comes the main lemma of this subsection,  crucial for the proof of the main
theorem of this article, Theorem \ref{Main}:

\begin{Lemma}\label{refer1}
Let $\{M_n\}_{n \in \N}$ and  $\{N_n\}_{n \in \N}$ be filtrations of the
CW-complexes $M$ and $N$, so that each $M_n$ (respectively $N_n$) is also a
CW-complex, included in $M$ (respectively $N$), cellularly, for each $n\in
\N$. Suppose  that $M$ is embedded in $N$, cellularly, and
that\footnote{Recall that if $A$ and $B$ are subcomplexes of a CW-complex $C$,
  then $A \cap B$ is a subcomplex of $C$ with a cell decomposition consisting
  of the cells of $C$ occurring in both $A$ and $B$.} $M_n=M \cap N_n,\forall n\in
\N$.   If, for any $n\in \N$, the inclusion of $M_n$ in $N_n$ is a homotopy
equivalence,   then there exists a deformation retraction from $N$  onto $ M$,
say $\r\colon N \times I \to N$, such that the restriction of $\r$ to $N_n
\times I$ is a deformation retraction of $N_n$ in $ M_n$, for any $ n \in
\N$. In fact $\r$ can be chosen to be cellular, in the sense that the map $\r\colon N\times I \to N$ is cellular.
\end{Lemma}
\begin{Definition} Deformation retractions with this property will be  called
  filtered. \end{Definition}
\begin{Proof} {\bf (Lemma \ref{refer1})}
Suppose first that the filtrations of $M$ and $N$ have length $L=1$. There
exists a deformation retraction $\r_0\colon N_0 \times I \to N_0$, of $N_0$ in
$M_0$, which we can suppose to be cellular. Let $r\colon N_1 \times I\to N_1$
be a deformation retraction of $N_1$  in $M_1$, also chosen to be  cellular. Consider the CW-complex $N_1 \times I$. Let us define a function $\r_1\colon N_1 \times I \to N_1$, in such a way that $\r_1$ is a cellular deformation retraction from $N_1$ onto $ M_1$, extending $\r_0$ (thus  permitting the argument to be extended to the general case by induction). We define $\r_1(m_1,t)=m_1, \forall m_1 \in M_1, \forall t \in I$. Define also $\r_1(n_0,t)=\r_0(n_0,t), \forall n_0 \in N_0, \forall t \in I $, and $\r_1(n_1,0)=n_1,\forall n_1 \in N_1$. There is no contradiction since $M_1 \cap N_0=M_0$, and $\r_0$ is a deformation retract of $N_0$ in $ M_0$.

Let us now extend $\r_1$ to all of $N_1 \times I$. Let $k \in \N$.
Suppose we have already extended the function $\r_1$ to  $(N_1^{k-1} \times I) \cup   (N_0 \times I) \cup (M_1 \times I) $, yielding a cellular map, with $\r_1(N_1^{k-1} \times \{1\} \cup   N_0 \times \{1\} \cup M_1\times \{1\})\subset M_1$.  Let $e^k$ be a $k$-cell of $N_1$ not belonging to $M_1 \cup N_0$. The function $\r_1$ is already defined in $\d e^k \times I$ and $e_k \times \{0\}$. Moreover $\r_1(\d e^k \times \{1\}) \subset M_1$. We can extend $\r_1$ to all $e_k \times I$, by using the cellular deformation retraction $r\colon  N_1 \times I \to N_1 $ (from $N_1$ to $M_1$). Let us be explicit. Choose a homeomorphism:
 \[f\colon \overline {e^k}\times  I \to (\d e^k \times I \cup e^k\times \{0\})\times I,\] 
such that $f(\overline{e^k} \times \{1\}) \subset  (\d e^k \times I \cup e^k \times \{0\})\times \{1\} \cup (\d e^k \times \{1\})\times I$, and $f$ sends  $(\d e^k \times I \cup e^k \times \{0\})\subset \overline{e^k} \times I$ identically to its copy $(\d e^k \times I \cup e^k\times \{0\}) \times \{0\} \subset (\d e^k \times I \cup e^k\times \{0\})\times I $. In particular $f$ is cellular. Then we extend $\r_1$ to $e^k\times I$ in the following way:

     \[\r_1(x,s)=r((\r_1 \times \id)\circ f(x,s)),x \in e^k,s \in I.\] 
Therefore $\r_1(e^k \times \{1\}) \subset M_1$, and the function obtained is still cellular. The result for $L=2$ follows from induction in $k$. As we have seen, the general case follows from the proof of the case $L=1$, by an inductive argument. 
\end{Proof} 

Let $M$ and $N$ be CW-complexes. Suppose that they have unique $0$-cells,
which we take to be their base points $*$. Give both $M$ and $N$  filtrations $\{M_n\}$ and $\{N_n\}$ by subcomplexes. Let $f\colon (M,*) \to (N,*)$ be a filtered map, which we suppose, furthermore, to be cellular. Consider the mapping cylinder:

\[\Mce_f=M \times I \bigcup_{(x,1) \mapsto f(x)} N,\] 
of $f$. Recall that $\Mce_f$ has a natural  deformation retraction onto $ N$,
say $\r\colon  \Mce_f \times I \to \Mce_f$, obtained by sliding the segment
$\{x\} \times I$  (where $x \in M$) along itself, towards the endpoint $f(x)
\in N$. Therefore, there exists a natural retraction $r\colon \Mce_f \to N$,
where $r(x)=\r(x,1), \forall x \in \Mce_f$; in other words $r(a,t)= f(a), \forall a \in M, \forall  t \in I$ and $r(b)=b, \forall b \in N$.
Note that  $f$ is $r$ composed with the inclusion map $i\colon M \to \Mce_f$.

 The mapping cylinder $\Mce_f$ of $f\colon M \to N$ has a natural cell decomposition, with $M$ and $N$ contained in $\Mce_f$, cellularly. Additionally, we have a $(n+1)$-cell of $\Mce_f$, for each $n$-cell $e^n$ of $N$, {\it connecting} $e^n$ with $f(e^n)$. This makes sense since $f$ is cellular.

The reduced mapping cylinder $\Mce_f'$ of $f$ is obtained from the CW-complex
$\Mce_f$ by collapsing the 1-cell of $\Mce_f$ connecting $*\in M$ with $*\in
N$ to a point. Therefore, $\Mce_f'$ has a unique $0$-cell (which we also call
$*$), and both $M$ and $N$ are included in $\Mce_f'$, cellularly, although
intersecting in $*\in \Mce_f'$. If $f$ is a homotopy equivalence, then $M$ is a deformation retract of $\Mce'_f$.

The CW-complex $\Mce'_f$ is filtered, in a natural way, by the reduced mapping cylinders $\Mce'_{f_n}$ of the restrictions $f_n\colon M_n \to N_n$ of $f$ to $M_n$, where $n \in \N$. Note that $\Mce'_{f_n}$ is included in $\Mce'_f$, cellularly, for any $n \in \N$. In addition,  both $M_n$ and $N_n$ are embedded in $\Mce'_{f_n}$, cellularly, and, moreover, $M_n=\Mce'_{f_n} \cap M$ and $N_n=\Mce'_{f_n} \cap N$, for any $n \in \N$.  From Lemma \ref{refer1} it follows that:

\begin{Corollary}\label{retraction}
Let $M$ and $N$ be CW-complexes with a unique 0-cell, provided with
filtrations  $\{M_n\}$ and $\{N_n\}$ by subcomplexes. Let $f\colon M \to N$ be a
filtered map,  which we assume to be  cellular, such that the restriction
$f_n$ of $f$ to $M_n$ is a homotopy equivalence $f_n\colon M_n \to N_n$.
There exist  filtered deformation retractions  of $\Mce'_f$ onto $M$
(respectively of 
$\Mce_f'$ onto $ N$), say $\r,\r'\colon \Mce'_f \times I\to \Mce'_f$ (respectively). Therefore the
restrictions of $\r$ and $\r'$ to $ \Mce'_{f_n} \times I$ are deformation
retractions from $\Mce'_{f_n}$ to $M_n$ (respectively  from $\Mce'_{f_n}$ to $N_n$), for any $n \in \N$. 
\end{Corollary}

Irrespectively of $f:M \to N$ being a homotopy equivalence or not, the
reduced mapping cylinder $\Mce_f'$ of $f$  always deformation retracts onto
$N$, in the same fashion  that the  mapping cylinder $\Mce_f$ of $f$ does. This
deformation retraction $\r\colon  \Mce_f' \times I \to \Mce_f'$, induced by the
natural  deformation retraction of  $\Mce_f$ onto $ N$,  is  filtered.
Therefore, the  obvious retraction $R\colon  \Mce_f' \to N$, such that
$R(x)=\r(x,1), \forall  x \in \Mce_f'$, thus $R(a,t)= f(a), \forall a \in M, \forall  t \in I$ and $R(b)=b, \forall b \in N$, is a filtered homotopy inverse of the
inclusion map $j:N \to \Mce_f'$. As before, we have $f=R \circ i$, where $i\colon  M
\to \Mce_f$ is the inclusion map. These very simple facts will be essential
later.

Note that Corollary \ref{retraction}  is valid without the assumption that the  complexes $M$ and $N$  have unique 0-cells, considering mapping cylinders instead of  reduced mapping cylinders. As a consequence, we obtain an extension of Theorem \ref{Simplify} for infinite filtrations, as long as they are filtrations by subcomplexes of CW-complexes.

\begin{Corollary}\label{Simplifyinfinite}
Let $M$ and $N$ be CW-complexes, equipped with filtrations  $\{M_n\}$ and $\{N_n\}$ by subcomplexes of $M$ and $N$, respectively. Let $f\colon M \to N$ be a filtered map, which we suppose, furthermore, to be cellular. Suppose that $f_n\colon M_n \to N_n$ is  a homotopy equivalence for any $n \in \N$. Then $f$ is a filtered homotopy equivalence.
\end{Corollary} 

\begin{Proof}
Consider the filtration of the mapping cylinder  $\Mce_f$ of $f$ by the mapping cylinders $\Mce_{f_n}$ of the restrictions  $f_n\colon M_n \to N_n$ of $f$ to $M_n$, where $ n \in \N$. Then, by the previous corollary, there exist filtered  homotopy equivalences $M \cong \Mce_f \cong N$. Specifically, the inclusions $i:M \to \Mce_f$ and $j:N \to \Mce_f$ are filtered homotopy equivalences. The obvious retraction $r: \Mce_f \to N$ is a filtered homotopy inverse of $j: N \to \Mce_f$. The result follows from the fact that $r \circ i=f$.
\end{Proof}

\subsection{Crossed Complexes}\label{ccomplex}
This subsection will only be needed in section $3$. We gather some results on
crossed complexes which we will use. Nothing here is new, and most is due to
R. Brown and  P.J. Higgins. The exceptions are   \ref{Tfree} and
\ref{Homotopy}, whose results are mainly  due to J.H.C. Whitehead.

Let $G$ and $E$ be groups. Recall that a  crossed module with base $G$ and fibre $E$, say $\G=(G,E,\d,\t)$, is given by a group morphism $\d \colon  E \to G$ and an action $\t$ of $G$ on $E$ on the left by automorphisms. The conditions on $\t$ and $\d$ are:
\begin{enumerate}
\item $\d(X \t e)=X\d(e)X^{-1},\forall X \in G, \forall e \in E$,
\item $\d(e) \t f=e f e^{-1}, \forall e, f \in E$.
\end{enumerate}
Notice that the second condition implies that $\ker \d$ commutes with all of  $E$.
We call $G$ the base group and $E$ the principal group. A morphism $F=(\f,\psi)$ between the crossed modules $\G$ and
$\G'=(G',E',\d',\t')$ is given by a pair of group morphisms $\f\colon G\to G'$ and
$\psi\colon E \to E'$, making the diagram
\begin{equation*}
\begin{CD}
E @>\p>> E' \\
@V\d VV  @VV\d' V\\
 G @>>\f > G'
\end{CD}
\end{equation*}
commutative. 
In addition we must have:   \[\f(X)\t' \psi (e)=\psi(X \t e), \forall X \in G, \forall e \in E.\]

There exists an extensive literature on crossed modules. We refer, for
example to 
\cite{BA1,BA2,B1,B4,BH1,BHu,BS,FM2}. A natural generalisation  of the concept of  a crossed module is a crossed complex:

\begin{Definition}
A (reduced) crossed complex $\A$ is given by a complex of groups:
\[...\to A_n \ra{\d_n=\d}  A_{n-1} \ra{\d_{n-1}=\d}A_{n-2} \to...\to A_2
\ra{\d_2=\d} A_1  \left (\ra{p} A \to \{1\}\right),\]
such that:
\begin{enumerate}

\item There exists a left action $\t=\t_n$ of the group $A_1$ on $A_n$, by automorphisms, for any $ n \in \N$, and all the boundary maps $\d$ are $A_1$-module morphisms. In addition, we suppose that $A_1$ acts on itself by conjugation.

\item The map  $A_2 \ra{\d_2} A_1$ together with the action $\t$ of $A_1$ in $A_2$ defines a crossed module. In other words to $1.$ we add the condition $\d(e) \t f=efe^{-1}, \forall e,f \in A_2$.
\item The group $A_n$ is abelian if $n>2$.

\item The action of $A_1$ on $A_n$  factors through the projection $p\colon  A_1 \to A=\coker(\d_2)$, for $n>2$, and, in particular, $A$ acts on $A_n,n>2$, on the left, by automorphisms.
\end{enumerate}
An $L$-truncated crossed complex is a crossed complex such that $A_n=\{0\}$ if $n>L$.
\end{Definition}

A natural example of a crossed complex is the following  one, introduced by
A.L. Blackers in \cite{Bl}:
\begin{Example}
Let $M$ be a path connected space, and let $\{M_n\}, n\in \N $ be a filtration
of it, where all spaces $M_n,n \in \N$ are path connected, and $M_0$ is a
singleton. Then the sequence of groups $\pi_n(M_n,M_{n-1},M_0=*)$, with the
obvious boundary maps, and left actions of $\pi_1(M_1,*)$ on them, is a crossed complex, which we denote by $\Pi(M)$, and call the ``Fundamental Crossed Complex of the Filtered Space $M$''. If the filtration of $M$ is finite, of length $L\in \N$, then, $\Pi(M)$ is an $L$-truncated crossed complex, which we denote by $\Pi_L(M,M_{L-1},...,M_2,M_1,*)$.
\end{Example}

\begin{Remark}
Let $M$ be a CW-complex with a unique $0$-cell, which we take to be  its base
point. The notation $\Pi(M)$ will always mean the fundamental crossed complex
of the skeletal filtration of $M$. If $L\in \N$, we will also denote $\Pi_L(M)
\doteq \Pi_L(M,M^{L-1},...,M^1,M^0=*)$.  This type of  crossed complexes was
considered  by J.H.C. Whitehead in \cite{W4,W5}.
\end{Remark}

It is easy to show  that crossed complexes and their morphisms, defined in the obvious way,
form a category. Crossed complexes are extensively studied or used in
\cite{BA1,BA2,B1,B4,BH2,BH3,BH4,BH5,BS,W4,W5}, for example. Notice that H.J. Baues
calls them ``crossed chain complexes''. On the other hand, Whitehead
J.H.C. only considered totally free crossed complexes, to be defined below,  referring to them as ``homotopy systems''.

We will usually denote a crossed complex $\A$ by $\A=(A_n,\d_n,\t_n)$, or more simply by $(A_n,\d_n)$, or  even $(A_n)$. A morphism $f\colon \A \to \B$ of crossed complexes will normally be denoted by $f=(f_n).$

 The category of crossed complexes  is a category with colimits. See \cite{BH2,BH3,BH4}. In particular we can consider the free product $\A \vee \B$ of two crossed complexes, defined as the pushout of the diagram:
\begin{equation*}
\begin{CD}
&\{1\} @>>> &\A\\
&@VVV &\\
&\B &
\end{CD},
\end{equation*}
and, analogously, for free products of infinite families of crossed complexes.
Here $\{1\}$ is the trivial crossed complex, so that $\{1\}_n=\{1\}, \forall n \in \N$. Therefore, if $\C$ is a crossed complex, there exists a one-to-one correspondence between $\Hom(\A \vee \B,\C)$ and $\Hom(\A,\C) \times \Hom(\B,\C)$.

 We will not need the  explicit construction of $\A \vee \B$, although we will
 use  the concept of free product  extensively.

\begin{Remark} \label{Generalise}
The definition of crossed complexes just given can be slightly generalised: Indeed, we 
can consider that $A_1$ is a groupoid with set of objects $C$, and each $A_n$
is a totally disconnected groupoid over the same set $C$. In addition all
group actions must be substituted by groupoid actions, in the obvious way. See \cite{B1,B3,BH2,BH3}, for
example, for a discussion of crossed complexes in this full generality, which
we will need to invoke later. We will use the designations ``reduced'' and ``non-reduced'', to distinguish between the group and groupoid-based definitions of crossed
complexes, respectively, whenever ambiguity may arise.

It is important to note that if $\A$ and $\B$ are reduced crossed complexes, then their free products in the category of reduced and non-reduced crossed complexes coincide. This also holds for general colimits.
\end{Remark} 

\subsubsection{Totally Free Crossed Complexes}\label{Tfree}

\begin{Definition}\label{freecrossedmodule}
Consider a group $G$. 
Let $K$ be a set and $\d_0\colon K \to G$ be a map.  A crossed module $\G=(G,E,\d,\t)$ is said to be the free crossed module  on $K$ and $\d_0\colon  K \to G$ if there exists an injective map $i\colon K \to E$ such that $\d \circ i=\d_0$, and the following universal property is satisfied:

\emph{
For any crossed module $\G'=(G',E',\d',\t')$, any  group morphism $\f\colon G \to G'$, and any map $\p_0\colon K \to E'$ for which $\d' \circ \p_0=\f \circ \d_0$, there exists a unique crossed module map $(\f,\p)\colon \G \to \G'$ such that $\p \circ i=\p_0$.}

\end{Definition}

For more details on the construction of these free crossed modules (defined up to isomorphism), we refer the reader, for example, to  \cite{BS,B4,BHu,FM2}. 

Recall the following theorem, due to J.H.C. Whitehead. For the original
proof see \cite{W1,W2,W4}.  See also \cite{BS}, $5.4$ and \cite{B3,BH2,GH}. 

\begin{Theorem} 
\label{Whitehead} {\bf (Whitehead's Theorem)}
Let $M$ be any path-connected  space with a base point $*$. Let $N$ be a  space obtained from $M$ by attaching some 2-cells  $s_1,...,s_n$.   Each 2-cell $s_i$ therefore induces an element $\f(s_i)$ of $\pi_1(M,*)$ (defined up to conjugation by an element of $\pi_1(M,*)$) through its attaching map. Then  the crossed module $\Pi_2(N,M,*)$ is the free  crossed module over the map $\f \colon \{s_1,...,s_n\} \to \pi_1(X,*)$.
\end{Theorem}

\begin{Definition}\label{totallyfree}
A (reduced) crossed complex $\A$ of the form: 
\[...\to A_n \ra{\d_n=\d}  A_{n-1} \ra{\d_{n-1}=\d}A_{n-2} \to...\to A_2 \ra{\d_2=\d} A_1 \ra{p} A \to \{1\},\]
is said to be totally free if:
\begin{enumerate}
\item $A_1$ is the free group on a set $C_1$.

\item The map $A_2 \ra{\d_2} A_1$, together with the action of $A_1$ on $A_2$,
  is the free crossed module on some map $\d^0_2\colon  C_2 \to A_1$. 

\item If $n\ge 3$ then  $A_n$ is the free $\Z(A)$-module on a set $C_n$. In particular  $\d_n$ is determined, in the obvious way, from a certain map $\d^0_n\colon C_n \to A_{n-1}$. 
\end{enumerate}
Notice that the action of an $a \in A_1$ on $\ker (\d_2)\subset A_2$ depends only on the projection $p(a)$ of $a$ in $\coker (\d_2)=A$, by the second condition of the definition of crossed complexes.

The sets $C_n, n\in \N$ will collectively be called a basis of $\A$.

\end{Definition}

Let $M$ be a CW-complex with a unique $0$-cell, which we take to be  its
base point. Then $\Pi(M)$ is a totally free crossed complex. This is proved, for
example, in \cite{W3,W4,B1}, and is a consequence of well known results.  Let us  be more explicit. For a more careful treatment of these issues we refer the reader to \cite{Br}. If $n \in \N$,  we can take each 
$C_n$ to be the set  $\{e^n_1,...,e^n_{l_n^M}\}$ of $n$-cells of $M$. 
 Each of these $n$-cells attaches to $M^{n-1}$ along an element $\d (e^n_{i_n}) \in
 \pi_{n-1}(M^{n-1},*), i_n=1,...,l_n^M$, defined up to acting by an element of
 $\pi_1(M^1,*)$. Therefore, if $n>2$, then  $\d^0_n$ is the image of this
 element on $\pi_{n-1}(M^{n-1},M^{n-2},*)$. It is well known that
 $\pi_n(M^n,M^{n-1},*)$ is the free $\Z(\pi_1(M^2,*))$-module on the set of $n$-cells
 of $M$. See \cite{GW}, chapter $V$, or \cite{B1}, for example.  We can see
 that, for dimensions higher than 3, a lot of information about the CW-complex
 $M$ is lost through passing to the fundamental crossed complex of its skeletal filtration.

If $n=2$, then $C_2$ is the set  of 2-cells of $M$, and, by Whitehead's Theorem, the crossed module $\Pi_2(M^2,M^1,*)$ is the free crossed module on the attaching maps $\d^0_2\colon C_2 \to \pi_1(M^1,*)$, defined up to conjugation by an element of $\pi_1(M^1,*)$.

We have:
\begin{Lemma} \label{colourings}
 Let $M$ be a CW-complex with a unique 0-cell which we take to be  its base
 point $*$. Let $C_n$ be the set of cells of $M$ of dimension $n$, where $n
 \in \N$. For for
 each $e^n \in C_n$, choose an element of $\pi_{n-1}(M^{n-1},*)$ (defined up
 to acting by an element of $\pi_1(M^1)$) along which $e^n$ attaches to $M^{n-1}$, therefore
 defining  maps $\d_n^0\colon C_n \to \pi_{n-1}(M^{n-1},M^{n-2},*), n \ge 3$ and
 $\d_2^0\colon  C_2 \to \pi_1(M^1,*)$. These choices also define identifications 
 $F(C_n) \to \pi_n(M^n,M^{n-1},*)$ and $F(C_1)\to \pi_1(M^1,*)$. Here $F(C_n)$
 is the free $\Z(\pi_1(M))$-module on $C_n$, if $n>3$, $F(C_2)$ is the
 principal group of the free crossed module on $\d_2^0\colon C_2 \to \pi_1(M^1,*)$,
 and $F(C_1)$ is the free group on $C_1$.

Let $\A=(A_n,\d_n)$ be a (reduced) crossed complex.  There exists a one-to-one correspondence between crossed complex morphisms $\Pi(M) \to \A$ and sequences of maps $f_n\colon  C_n \to A_n$ such that $f_{n-1} \circ \d^n_0=\d_n \circ f_n$ for any $n \in \N$. 
\end{Lemma}

\subsubsection{The Cotruncation Functor} 
Let $M$ be a CW-complex.  As usual,  for any $k \in \N$, let $M^k$ denote the
$k$-skeleton of $M$. Suppose that $M$ has a unique $0$-cell, which we take to
be its base point $*$.  Let $L>1$. Note that from  the Cellular Approximation Theorem we have: 
 \[\Pi_{L}(M,M^{L-1},M^{L-2},...,*) \cong \Pi_L(M^{L+1},M^{L-1},M^{L-2},...,*) .\]

Consider the fundamental crossed complex: 
\[... \ra{\d_4} \pi_3 (M^3,M^2,*) \ra{\d_3} \pi_2(M^2,M^1,*) \ra{\d_2}
\pi_1(M^1,*)\ra{p} \pi_1(M,*) \to \{1\},\]
of the skeletal filtration of $M$.
 The following is well known:
\begin{Lemma} Let $L>1$.
We have:
\[\pi_L(M,M^{L-1},*)=\pi_L(M^L,M^{L-1},*)/{\rm im}(\d_{L+1}).\]
\end{Lemma}
This Lemma can also be shown using the General van Kampen Theorem (see Theorem \ref{vkt}).

\begin{Proof}
This follows from the homotopy sequence of $(M^{L+1},M^{L},M^{L-1})$. Indeed the following sequence is exact:
\begin{multline*}
  \to \pi_{L+1}(M^{L+1},M^L,*) \ra{\d_{L+1}} \pi_L (M^{L},M^{L-1},*) \to \pi_L(M^{L+1},M^{L-1},*)\\ \to \pi_L(M^{L+1},M^L,*)\to ,
\end{multline*} 
and, moreover, $\pi_L(M^{L+1},M^{L},*) \cong \{1\}$.  Recall that  $\pi_L(M,M^{L-1},*)=\pi_L(M^{L+1},M^{L-1},*)$.
\end{Proof}

For any $L \in \N$, there exists a cotruncation functor $\CT_L$ from the
category of crossed complexes to the category of $L$-truncated crossed
complexes. If $\A=(A_n,\d_n,\t_n)$ is a crossed complex, then $\CT_L(\A)$ is equal to $\A$ up to dimension $L-1$, whereas $A_L$ is substituted by $A_L/{\rm im}(\d_{L+1})$. Therefore, we can rephrase our previous lemma as:
\begin{Proposition}\label{Cotrunc}
Let $M$ be a CW-complex with a unique $0$-cell, which we take to be its base point. Let $L \in \N$. We have: 

\[\Pi_L(M,M^{L-1},M^{L-2},...,*) = \CT_L(\Pi(M)).                                 \]
\end{Proposition} 

\subsubsection{Colimit Theorems: The ``General van Kampen Theorem''}
Let $M$ be a connected space equipped with a filtration $\{M_n\}_{n \in \N}$.
\begin{Definition}
We say that $\{M_n\}_{n \in \N}$ is connected if all the spaces $M_n, n \in
\N$ are connected, and, moreover, $\pi_n(M_k,M_n,*)=0$ if $k\ge
n$, for any  $n \in \N$. 
\end{Definition}
Therefore, if $M$ is a connected CW-complex with a unique 0-cell, then its
skeletal filtration is connected. Similarly, if $L>1$, then any filtration of the
type $(M,M^L,M^{L-1},...,M^0)$ is connected as well.
The following beautiful theorem is due to R. Brown and P.J. Higgins See
\cite{BH2,BH3,BS,B1,B4}. It is called ``General van Kampen Theorem''.
\begin{Theorem} {\bf (R. Brown and P.J. Higgins)}\label{vkt}
Let $M=\{M_n\}$ and  $N=\{N_n\}$  be filtered spaces. Let $U\subset M$, and let $i\colon U \to M$ be
the inclusion map. Let $U$ inherit the filtration $\{U_n\}$
induced by the filtration of $M$. Suppose that the filtrations of $M,N$ and
$U$ are connected. In addition, we suppose that the inclusions $i_n\doteq
i_{|U_n} \colon U_n \to M_n$
are closed cofibrations, for any $n \in \N$. Let $f\colon U \subset M \to N$ be a filtered map. Therefore the
adjunction space:
 \[V\doteq M \bigcup_{u \in U \mapsto f(u)} N\] 
is naturally filtered
by the sets:
  \[V_n\doteq M_n \bigcup_{u \in U_n \mapsto f(u)} N_n, n \in \N.\]
We thus have a commutative diagram of filtered spaces:
\begin{equation*}
\begin{CD}
&U @>f>> &N\\
& @V iVV &@VVi'V \quad .\\
&M @>f'>> &V\\
\end{CD}
\end{equation*}
The following holds:

\begin{enumerate}
\item The filtered space $V$ is connected.
\item The diagram:
\begin{equation*}
\begin{CD}
&\Pi(U) @>f_*>> & \Pi(N)\\
& @Vi_*VV &@VVi'_*V\\
&\Pi(M) @>>f'_*> &\Pi(V)\\
\end{CD}
\end{equation*}
\end{enumerate}
is a pushout of crossed complexes.
\end{Theorem}

This theorem also holds for  infinite adjunction spaces. In particular it follows:
\begin{Corollary}\label{Free}
Let $\{M_\l, \l \in \L\}$ be a family of CW-complexes. Suppose that they all
have unique 0-cells, which we take to be their base points $*$. 
Form the wedge  product:
\[\V_{\l \in \L} (M_\l,*),\]
 along $*$, which is also a CW-complex $\V_{\l \in \L}  M_\l$, having a unique 0-cell.  The obvious map: 
\[\V_{\l \in \L} \Pi(M_\l)\to \Pi\left (\V_{\l \in \L} M_\l\right ),\] induced by the inclusions,  is an isomorphism of crossed complexes.  In addition,  if $L \in \N$, the natural map:  
\[\V_{\l \in \L} \Pi_L(M_\l) \to \Pi_L\left (\V_{\l \in \L} M_\l\right),\] is also an isomorphism.
\end{Corollary}
This is Theorem $7.1$ of \cite{BH2}.
\begin{Remark} The fact that the fundamental crossed complex of a CW-complex
  is totally free (referred in \ref{Tfree}) is also a corollary of the General
  van Kampen Theorem.
\end{Remark} 
\subsubsection{Homotopy of Crossed Complexes}\label{Homotopy}

The main references now  are \cite{B1,BH3,BH4,BH5,W4,W5}.
Let $G$ and $G'$ be groups. Suppose that  $G'$ has a left  action $\t'$ on the
group $E'$ by automorphisms. Let $\f\colon G \to G'$ be a group
morphism. A map $s\colon G \to E'$ is said to be a $\f$-derivation if:
\[s(XY)=\left (\f(Y)^{-1} \t' s(X)\right) s(Y), \forall X,Y \in G.\] 
Notice that, if $G$ is a free group, then a $\f$-derivation can be specified, uniquely,  by
its value on the set of free generators of $G$. This is easy to show, but see \cite{W4}, Lemma $3$, for example.

\begin{Definition}
Let $\A=(A_n,\d_n, \t_n)$ and $\B=(B_n,\d'_n,\t'_n)$ be crossed complexes, and let
$f=(f_n)\colon  \A \to \B$ be a morphism. An (L-fold) $f$-homotopy ($L\ge 1$) is a sequence $H=(H_n)$ of maps $H_n\colon 
A_n \to B_{n+L}$ such that:
\begin{enumerate}
\item If $n>1$ then $H_n$ is a group morphism  and, moreover:
\[H_n(a_1 \t a_n)=f_1(a_1) \t'H_n(a_n), \forall a_1 \in A_1, \forall a_n \in
A_n.\]
\item $H_1$ is an  $f_1$-derivation. 
\end{enumerate}
A (1-fold) $f$-homotopy will be simply called  an $f$-homotopy.
Note that there are no compatibility relations between the boundary maps  of $\A$ and $\B$ and
the maps $H_n$, where  $n \in \N$, whenever $H=(H_n)$ is an $L$-fold $f$-homotopy.
\end{Definition}
We have (see \cite{W4,B1,BH4}):
\begin{Lemma}\label{connection}
Under the conditions of the previous definition, if $H$ is an  $f$-homotopy, where $f\colon  \A \to \B$, then the sequence of maps
$g=(g_n)$ such that:
\[g_n(a_n)=f_n(a_n)(H_{n-1}\circ \d)(a_n) (\d' \circ H_n)(a_n), a_n \in A_n, n\ge
2,\] 
and  
\[g_1(a_1)=f_1(a_1)(\d' \circ H_1)(a_1), a_1 \in A_1,\]
is a morphism of crossed complexes, in which case, we say that $H$ realises a homotopy $f \ra{H} g$. 
\end{Lemma}

We can consider the category of crossed complexes with up to
homotopy 
morphisms of crossed complexes as morphisms. The following result is due to J.H.C. Whitehead, see \cite{W4} (Theorem $5$), 
and gives us one answer as to  how $\Pi(M)$ depends on the CW-decomposition of
$M$, if $M$ is a CW-complex.
\begin{Theorem}\label{homotopy} {\bf (J.H.C. Whitehead)}
Let $M$ and $N$ be CW-complexes with a unique 0-cell. Let $F,G\colon M \to N$ be homotopic cellular maps. Then the induced maps $F_*,G_*\colon  \Pi(M) \to \Pi(N)$ are homotopic. In particular, If $M$ and $N$ are
homotopic (as spaces) then it follows that $\Pi(M)$ and $\Pi(N)$ have the same homotopy type.
\end{Theorem}
This result can be proved using R. Brown and P.J. Higgins framework on $\w$-groupoids and the tensor product of crossed complexes, and is, ultimately, a consequence of the Cellular Approximation Theorem.  See \cite{BH4,BH5}. In fact, it is not difficult to show that a homotopy $H:F \to G$, chosen to be cellular, will induce a crossed complex homotopy $F_* \to G_*$.
\begin{Remark}
In fact, Whitehead's results imply a stronger statement:  Suppose  that  $M$ and $N$ are finite of dimension $\leq L$. Then   $\Pi(M')$ and $\Pi(N')$ are simply homotopy equivalent. Here $M'$ and $N'$ are obtained from $M$ and $N$ by making a wedge product with a certain finite number of spheres $S^L$. This is a consequence of Theorem $16$ of \cite{W5}. For the definition of simple equivalence of crossed complexes, we refer the reader  to \cite{W5}. 
\end{Remark}

One of our main results (Theorem \ref{Main2}) extends (J.H.C. Whitehead's) Theorem \ref{homotopy} in a direction very similar to the one indicated in the previous remark.
\subsubsection{The Definition of $CRS(\A,\B)$ }

Given crossed complexes $\A$ and $\B$,
there exists a groupoid $CRS_1(\A,\B)$ whose objects are the morphisms of
crossed complexes $\A \to \B$, and whose morphisms are the homotopies connecting them, in
the manner  shown in Lemma \ref{connection}, with pointwise product of maps as composition.  We
denote the set of all $f$-homotopies  by $CRS_1^f(\A,\B)$. Therefore
$CRS^f_1(\A,\B)$ is the set of morphisms of $CRS_1(\A,\B)$ with source $f\in \Hom(\A,\B)$.

In fact, a  far  stronger result is true: Let  $\A$ and $\B$ be crossed
complexes, and let  $f=(f_n)\colon \A \to \B$ be a morphism. Let $n\ge 2$. An $n$-fold
$f$-homotopy $H$ determines naturally a $(n-1)$-fold $f$-homotopy $\d^f_n(H)$,
in the same way that a 1-fold $f$-homotopy determines a morphism. Let
$CRS_n^f(\A,\B)$ be the set of $n$-fold $f$-homotopies. It is  a group with
pointwise multiplication as product. In particular for any $n>1$:
  \[CRS_n(A,B)=\bigcup_{f \in
  \Hom(\A,\B)} CRS^f_n(\A,\B)\] is a totally disconnected groupoid with the same
object set as $CRS_1(\A,\B)$. There exist also natural left groupoid actions of
$CRS_1(\A,\B)$ on $CRS_n(\A,\B)$, for any $n \ge 2$.  The following theorem is shown in \cite{B1,BH3,BH4,BH5}, where the details of
this description can be found.

\begin{Theorem} {\bf (R. Brown and P.J. Higgins)}
Let $\A$ and $\B$ be (reduced) crossed complexes. Then the sequence
$CRS(\A,\B)=(CRS_n(\A,\B))$ is a (non reduced) crossed complex.
\end{Theorem}

\begin{Remark}
Actually R. Brown and P.J. Higgins considered a more general notion of homotopy. Our more
restricted notion was called pointed homotopy by them. We will not need to
consider the full generality.
\end{Remark}

\subsubsection{Counting Homotopies} \label{conthom}
The following lemma will be extremely useful. It is due to the fact that there are
no constraint relations between a homotopy and the boundary maps of the crossed
complexes involved. It appears in \cite{BH4}.
\begin{Lemma}\label{extension1}
Let $\Fc=(F_n,\d_n)$ be a totally free (reduced) crossed complex. Let $C_n$ denote the
set of free generators of $F_n,n=1,2,...$, as a free group if $n=1$, as a free
crossed module if $n=2$, and as a free $\Z(\coker(\d_2))$-module if $n>2$. Let
$\A=(A_n,\d'_n)$ be a (reduced) crossed complex, and let $f\colon \Fc \to\A$ be a crossed
complex morphism. An
$L$-fold $f$-homotopy can be specified, uniquely, by its value on the sets
$C_n, n \in \N$.
\end{Lemma}  
\begin{Proof}
The only bit that needs discussion is the definition of $H_2\colon F_2 \to A_{L+2}$,
from  the known form of $H_2$ restricted to ${C_2}$. This follows from the
universal property defining  free  crossed modules (Definition \ref{freecrossedmodule}), applied to the crossed module $A_{2+L}\ra{a \mapsto 1} A={\rm coker} (\d'_2)$. Note condition $4.$ of the definition of crossed complexes. 
\end{Proof}

Therefore:

\begin{Corollary}\label{extension}
 Let $M$ be a CW-complex with a unique $0$-cell, which we take to be its base point. Let also $\A$ be a reduced crossed complex, and let $f\colon \Pi(M)\to \A$ be a
  morphism of crossed complexes. An $n$-fold f-homotopy is uniquely specified
  by its value on each cell of $M$.  In fact, if $\A$ is $L$-truncated, then
  an $n$-fold homotopy $\Pi_L(M)$ to $M$ is also uniquely specified by its value on the cells of $M$ of dimension $1,2,...,L-1$. 
\end{Corollary}

The following simple quantification of the previous corollary will have an important role later:

\begin{Corollary}\label{cardinality}
 Let $M$ be a  CW-complex with a unique $0$-cell, which we take to be its base point, and a finite number of $n$-cells for any $n \in \N$. 
Let also $\A=(A_n)$ be a (reduced) crossed complex, which we
  suppose to be finite, in the sense that $\A$ is $l$-truncated, for some
  $l\in \N$, and each $A_n,n \in \N$ is finite. Then $CRS(\Pi(M),\A)$ is a
  finite crossed complex  and, moreover, for any $L \in \N$ and any $f \in
  \Hom\left (\Pi(M),\A\right )$ we have:
\[\#(CRS^f_L(\Pi(M),\A))=\prod_{k=1}^{\infty}\left ( \# (A_{k+L}) \right)^{l^M_k}, \]
where $l^M_k$ denotes the number of cells of $M$ of order $k\in \N$.
\end{Corollary}

\subsubsection{The Classifying Space of a Crossed Complex}\label{Classifying}
There exists a functor, the classifying space functor, from the category of
(reduced or non reduced) crossed complexes to the category of filtered spaces. This functor, due to
R. Brown and P.J. Higgins, was described in \cite{BH5}.   It generalises the
concept of the classifying space of a group. A similar construction of classifying spaces of crossed complexes appeared in \cite{Bl}.

For a crossed complex $\A=(A_n)$, let $|\A|$ denote its classifying space. It is a
CW-complex, with one 0-cell for each element of $C$, where $C$ is the set of
objects of the groupoid $A_1$. A very strong  result appearing in  \cite{BH5} is the following one, which is a consequence of well know simplicial techniques:

\begin{Theorem} {\bf (R. Brown and P.J. Higgins)}\label{weaktype}
Let $M$ be a CW-complex, provided  with its skeletal filtration. Suppose $M$ has a unique 0-cell which we take to be its base point.  Let also $\A$ be a
 reduced crossed complex, thus $|\A|$ has a unique $0$-cell $*$. There exists a map: \[\p\colon |{CRS(\Pi(M),\A)}| \ra[\simeq_w]{} TOP((M,*),(|\A|,*)),\]
 which is  a weak homotopy equivalence.
\end{Theorem} 
In particular, if $M$ is a finite CW-complex then $\p$ is a homotopy equivalence. This is because, in this case,  the function space $TOP((M,*),(|\A|,*))$ has the homotopy type of a CW-complex. See \cite{Mi}.

Let $\A$ be a crossed  complex. We will consider now the full generality in
the definition of crossed complexes, as sketched in Remark
\ref{Generalise}. Let $C$ be the object set of the groupoid $A_1$. Denote the set of morphisms of $A_n$ whose source is $c$  by
$A^c_n$ . Therefore  $A^c_n$  is a group
if $n>1$. Let also $A_1^{(c,d)}$ be the set of morphisms of $A_1$ with source
and target $c,d \in C$, respectively. Let $\d^c_n={\d_n}_{|A^c_n}$. Consider the  complex $\A^c=(A_n^c,
\d^c_n),n>1$, considering the group $\A_1^{c,c}$ at index $1$.  Define
$\pi_1(\A,c)=H_1(\A^c)$, and $H_n(\A,c)=H_n(\A^c)$ if $n\ge 2$.  Notice that
these homotopy and homology groups with base $c\in C$ only depend on the
connected component in the groupoid $A_1$ to which $c$ belongs.

If $M$ is a CW-complex, with a unique 0-cell, then $\pi_1(M,*)\cong \pi_1(\Pi(M),*)$, whereas $H_n(\M) \cong H_n(\Pi(M))$, for $n \ge 2$. Here $\M$ is the universal covering of $M$.

\begin{Theorem}\label{homgroups}{\bf (R. Brown  and P.J. Higgins)}
Let $\A$ be a crossed complex. Let $C$ be the set of objects of the groupoid $A_1$. Then  the classifying space $|\A|$ of $\A$ is a CW-complex with one zero cell $c$ for each element of $c \in C$. Moreover:
\begin{align*}
\pi_1(|\A|,c)&\cong\pi_1(\A,c)\\
\pi_n(|\A|,c)&\cong H_n(\A,c),n\ge 2,\\
\end{align*}
for any $c \in C$. Furthermore, there exists a one-to-one correspondence between
connected components of the groupoid $A_1$ and connected components of $|\A|$.

\end{Theorem}

\section{The Homotopy Type of the Skeletal Filtration of a CW-Complex} 

Let $M$ be a space  which can be given a  CW-complex structure. In this
chapter we analyse to which extend the homotopy type of the skeletal filtration $\{ M^n\}_{n \in \N}$ of $M$ depends on the cellular decomposition of $M$. 

Throughout all this chapter, $D^n$ denotes the $n$-disk and $S^n$ denotes the $n$-sphere. As usual we set $I=[0,1]$.

\subsection{Dimension Two}\label{two} 
Let $(N,M)$  be a pair of CW-complexes  such that the inclusion of $M$ in $N$
is a homotopy equivalence. For simplicity, assume that  $M$ and $N$ have a
finite number of 1-cells. Our whole discussion  will remain valid without this restriction, with the obvious modifications. 

 Let $M^1$ and $N^1$ be, respectively, the 1-skeletons of $M$ and $N$. Suppose that  $N$ and  $M$ have a unique 0-cell, which we take to be their common base point $*$, so that both $M$ and $N$ are well pointed. 

The group $\pi_1(M^1,*)$ is the free group on the set $\{X_1,...,X_m\}$ of
1-cells of $M$.  Let   $Y_1,...,Y_n$ be the 1-cells of $N$ which are not in $M$. Then $\pi_1(N^1,*)$ is the free group $F(X_1,...,X_m,Y_1,...,Y_n)$ on the set    $\{X_1,...,X_m,Y_1,...,Y_n\}$.

\begin{Theorem}
There exists a homotopy equivalence:
\[(N,N^1,*)\cong (M,M^1,*) \vee (D^2,S^1,*)^{\vee n}.\]                   
 \end{Theorem}

\begin{Proof} \footnote{This argument arose in a discussion with Gustavo Granja.}
Since $M$ is a subcomplex of $N$, and $N$ is homotopic to $M$, it follows that $M$ is a strong deformation retract of $N$. By the Cellular Approximation Theorem, we can suppose, furthermore, that there exists a retraction $r\colon N \to M$ sending $N^1$ to $ M^1$, and such that $r \cong \id_N$, relative to $M$.  In particular if $k \in \{1,...,n\}$ then we have $Y_k r_*(Y_k)^{-1}=1_{\pi_1 (N,*)}$, (though this relation does not  hold in $\pi_1(N^1,*)$).
Define a map \[f\colon (P,P^1,*) \doteq (M,M^1,*) \vee \bigvee_{k=1}^n (D^2_k,S^1_k,*)\to (N,N^1,*)\] in the following way: First of all, send $(M,M^1,*)$ identically to its copy $(M,M^1,*)\subset (N,N^1,*)$. Then we can send each $(S^1_k,*),k=1,...,n$ to the element $Y_k r_*(Y_k)^{-1}\in \pi_1(N^1,*)$. Since these elements are null homotopic in $(N,*)$, this map extends to the remaining 2-cells of $(P,P^1,*)$.

 Let us prove that $f\colon (P,P^1,*) \to (N,N^1,*)$ is a homotopy equivalence. It suffices to prove that $f\colon (P,*) \to (N,*)$ and $f^1\doteq f_{|P^1}\colon (P^1,*) \to (N^1,*)$ are based homotopy equivalences, by Theorem \ref{Simplify}. 
 
 We first show that $f$ is a homotopy equivalence $(P,*) \to (N,*)$.
Let $r'\colon (P,*)\to (M,*)$ be the obvious retraction, thus $r' \cong \id_{(P,*)}$. We have $r \circ f\cong r \circ f \circ r' =r'\cong \id_{(P,*)}$. On the other hand $f \circ r=r \cong \id_{(N,*)}$.

We  now show that $f^1$ is a  homotopy equivalence $(P^1,*) \to (N^1,*)$. It is enough to prove that the induced map $f^1_*\colon  \pi_1(P^1,*) \to \pi_1(N^1,*)$ is an isomorphism. Note that  $\pi_1(P^1,*)$ is (similarly with $\pi_1(N^1,*)$) isomorphic  to the free group on the set $\{X_1,...,X_m,Y_1,...,Y_n\}$. The induced map on the fundamental groups has the form: 
\[f^1_*(X_k)=X_k, k=1,...,m, \textrm{ and }f^1_*(Y_k)=Y_kr_*(Y_k)^{-1}, k=1,...,n.\]
 Notice that $r_*(Y_k) \in F(X_1,...,X_m), k=1,...,n$. Consider the morphism $g$ of $F( X_1,...,X_m,Y_1,...,Y_n)$ on itself such that:
\[g(X_k)=X_k, k=1,...,m \textrm{ and } g(Y_k)=Y_kr_*(Y_k), k=1,...,n.\]
 Therefore $(f^1_* \circ g )(X_k)=X_k, k=1,...,m,$ and:
\begin{align*}
 (f^1_* \circ g)(Y_k)&=f^1_*(Y_k r_*(Y_k))\\
                   &=f^1_*(Y_k) f^1_*(r_*(Y_k))\\
                   &=Y_kr_*(Y_k)^{-1}r_*(Y_k)\\
                   &=Y_k,k=1,...,n.\\
\end{align*}
 Analogously $(g \circ f_*)(X_k)=X_k, k=1,...,m$, and: 
\begin{align*}
 (g \circ f^1_*)(Y_k)&=g(Y_k r_*(Y_k)^{-1})\\
                     &=Y_k r_*(Y_k) g(r_*(Y_k^{-1}))\\
                     &= Y_k r_*(Y_k) r_*(Y_k^{-1})=Y_k,k=1,...,n.
\end{align*}
 This  proves that $g^{-1}=f^1_*$, which finishes the proof.
\end{Proof}

\begin{Corollary}
Let $M$ and $N$ be  CW-complexes, with  unique $0$-cells, which we take to be their base points, both denoted by $*$. Let $m_1$ and $n_1$ be the number of $1$-cells of $M$ and $N$, respectively. Suppose that  $M$ and $N$ are based homotopic, as spaces.  There exists a (filtered) homotopy equivalence:
\[f\colon  (M,M^1,*)\vee (D^2,S^1,*)^{\vee n_1} \to   (N,N^1,*)\vee (D^2,S^1,*)^{\vee m_1}.\]
Moreover, $f$ can be taken to be cellular.
\end{Corollary}
\begin{Proof}
Let $F\colon M \to N$ be a pointed homotopy equivalence. We can suppose that  $F$ is a cellular map. The reduced mapping cylinder $P$ of $F$ is a CW-complex with a unique 0-cell, containing $M$ and $N$ as subcomplexes. Moreover $P$ is homotopic to both $M$ and $N$. The complex $P$ has $n_1$ 1-cells which are not in $M$ and $m_1$ 1-cells which are not in $N$. By the previous theorem, it thus follows that there exists a filtered homotopy equivalence:
\[F'\colon \Mc \doteq (M,M^1,*)\vee (D^2,S^1,*)^{\vee n_1}\to  (N,N^1,*)\vee (D^2,S^1,*')^{\vee m_1} \doteq \Nc.\]

Let $f$ be a map, homotopic to $F'$,  which is cellular. We can suppose that  $F$ coincides with $F'$ in $(M^1,*)\vee (S^1,*)^{\vee n_1}=(\Mc^1,*)$. Then $f$ is a homotopy equivalence $f\colon \Mc \to \Nc$, and, therefore, by Theorem \ref{Simplify}, $f$ is a filtered homotopy equivalence.
\end{Proof}

\subsection{The General Case}\label{general}

Let $M$ and $N$ be CW-complexes, homotopic as (based)  spaces. For simplicity,
we suppose that they only have a finite number of cells at each dimension. The
main results that we obtain in this subsection (theorems \ref{Main},
\ref{MainInfty} and \ref{claim1}) are valid without this restriction, although they need to be
re-written in the natural  way. The proof we will make will remain valid, as
long as we make the obvious modifications. As usual, we suppose $M$ and $N$ to have a
unique 0-cell, which we take to be  the  base point of both,  so that  $M$ and $N$ are well pointed.

For a $k \in \N$, let $l^M_k$ and $l^N_k$ be the number of cells of $M$ and $N$ (respectively) of dimension $k$.
Define the integer numbers:
\begin{align*}
l^M_{\{M,N\}}(1)&=l^N_1,\\
l^N_{\{M,N\}}(1)&=l^M_1.
\end{align*}
And, in general, if $k>1$: 
\begin{align*}
 l^M_{\{M,N\}}(k+1)&=l^N_{\{M,N\}}(k)+l^N_{k+1},\\
 l^N_{\{M,N\}}(k+1)&=l^M_{\{M,N\}}(k)+l^M_{k+1}.
\end{align*}

Define, for each $L \ge 2$, the filtered space:
\begin{multline*}
\Mc_L=(M,M^{L-1},...,M^1,*) \vee \V_{j_{(L-1)}=1}^{l^M_{\{M,N\}}(L-1)}(D^L_{j_{(L-1)}},S^{(L-1)}_{j_{(L-1)}},*,...,*)\vee... 
\\... \vee \V_{j_{(L-2)}=1}^{l^M_{\{M,N\}}(L-2)}(D^{(L-1)}_{j_{(L-2)}},D^{(L-1)}_{j_{(L-2)}},S^{(L-2)}_{j_{(L-2)}},*,...,*)\vee... 
\\... \vee \V_{j_1=1}^{l^M_{\{M,N\}}(1)} (D^2_{j_1},...,D^2_{j_1},S^1_ {j_1},*).
\end{multline*}
Analogously, define:
\begin{multline*}
\Nc_L=(N,N^{L-1},...,N^1,*) \vee \V_{i_{(L-1)}=1}^{  l^N_{\{M,N\}}(L-1)   }(D^L_{i_{(L-1)}},S^{(L-1)}_{i_{(L-1)}},*,...,*)\vee... \\ ...\vee \V_{i_{(L-2)}=1}^{l^N_{\{M,N\}}(L-2)}   (D^{(L-1)}_{i_{(L-2)}},D^{(L-1)}_{i_{(L-2)}},S^{(L-2)}_{i_{(L-2)}},*,...,*)\vee... \\... \vee \V_{i_1=1}^{l^N_{\{M,N\}}(1)} (D^2_{i_1},...,D^2_{i_1},S^1_{i_1},*).
\end{multline*}
Therefore $\Mc_L$ has $l^M_L+l^M_{\{M,N\}}(L-1)=l^N_{\{M,N\}}(L)$ cells of dimension $L$, and, analogously,  $\Nc_L$ has $l^N_L+l^N_{\{M,N\}}(L-1)=l^M_{\{M,N\}}(L)$ cells of dimension $L$, for $L=2,3,...$.
Define also $\Mc_1=(M,*)$ and $\Nc_1=(N,*).$

For any $L \in \N$, the CW-complex $\Mc_L$ is embedded cellularly in $\Mc_{L+1}$,  and analogously for $\Nc_L$. Therefore the  spaces:
\begin{align*}
\Mc&=\bigcup_{L \in \N} \Mc_L\\
   &=M \vee \V_{L=1}^\infty \left (\V_{i_L=1}^{l^M_{\{M,N\}}(L)} D^{L+1}_{i_L}\right),
\end{align*}
and:
\begin{align*} 
\Nc&=\bigcup_{L \in \N} \Mc_L\\
   &=N \vee \V_{L=1}^\infty \left (\V_{j_L=1}^{l^N_{\{M,N\}}(L)}D^{L+1}_{j_L}\right)
\end{align*}
are CW-complexes.  Filter $\Mc$ and $\Nc$ by their skeletal filtrations.
Notice that we regard each  $(L+1)$-disk $D^{L+1}$, where $ L \in \N$, as having the obvious  CW-decomposition with a unique $0$-cell, an $L$-cell and an $(L+1)$-cell.

 It is important to note that the following holds for any $L \in \N$:
\begin{align*}
\Mc^{L-1}&=\Mc_L^{L-1},\\
\Nc^{L-1}&=\Nc_L^{L-1}.
\end{align*}

We want to prove the following theorem:
\begin{Theorem}\label{Main}
In the above situation,
for any $L\in \N$ there exists a (filtered) homotopy equivalence:
\[F_{L}\colon \Nc_L \ra[\cong]{} \Mc_L,\]
which can be taken cellular. Moreover, if we are provided with a filtered homotopy equivalence $F_{L}\colon \Mc_{L} \to \Nc_{L}$, then there exists a filtered homotopy equivalence $F_{L+1}\colon  \Mc_{L+1} \to \Nc_{L+1}$ which extends the map $F_{L}^{L-1}$, the restriction of ${F_{L}}$ to ${\Mc_{L}^{L-1}}=\Mc_{L+1}^{L-1}$. 
\end{Theorem}  

Since $\Mc_{L}^{L-1}=\Mc^{L-1}$ and $\Nc_{L}^{L-1}=\Nc^{L-1}$ for any $ L \in \N$, it follows (using Corollary \ref{Simplifyinfinite}):
\begin{Theorem}\label{MainInfty}
Let $M$ and $N$ be CW-complexes with unique 0-cells, taken to be their
bases points. Suppose that $M$ and $N$ are homotopic as (based) spaces. There exists a filtered homotopy equivalence:
\[F:\Mc \ra[\cong]{} \Nc.\]
\end{Theorem} 

The proof of  Theorem \ref{Main} is done by an obvious induction on $L$. Notice that we have already proved it for $L=2$, dealt with in \ref{two}. The general proof follows immediately from Theorem \ref{claim1} below,  and the already stated fact that $\Mc_L$ has $l^M_L+l^M_{\{M,N\}}(L-1)=l^N_{\{M,N\}}(L)$ cells of dimension $L$, for any $L \in \N$,  and analogously for $\Nc_L$.

\subsubsection{An Auxiliary Discussion}
Let $M$ and $N$ be CW-complexes. As usual, we suppose that they have unique
0-cells, which we take to be  their base points, both denoted by $*$. Suppose that for some $L> 1$ there exists a (filtered) homotopy equivalence:
\[F\colon (M,M^{L-1},M^{L-2},...,M^0=*) \to (N,N^{L-1},N^{L-2},...,N^0=*),\]
 which, furthermore,  is a cellular map. Let $m_L$ and $n_L$ be, respectively, the number of cells of $M$ and $N$  of dimension $L$.
\begin{Theorem}\label{claim1}
There exists a (filtered) homotopy equivalence:
\begin{multline*}
F''\colon (M,M^L,M^{L-1},....,M^0=*) \vee \V_{i_L=1}^{n_L}
(D^{L+1}_{i_L},S^L_{i_L},*,..,*)\doteq \bar{M} \\
\to (N,N^L,N^{L-1},....,N^0=*) \vee \V_{j_L=1}^{m_L}
(D^{L+1}_{j_L},S^L_{j_L},*,..,*)\doteq \bar{N},
\end{multline*}
which agrees with $F$ over $M^{L-1}$. Moreover, $F''$ can be taken cellular.
\end{Theorem}

\begin{Remark}
  This result, as well as Theorem \ref{Main},  is reminiscent of the following result of J.H.C. Whitehead on simple homotopy types: If $A$ and $B$ are finite CW-complexes of dimension $n$ with the same $(n-1)$-type, and if $f\colon A \to B$ realises an $(n-1)$-homotopy equivalence, then there exists a simple homotopy equivalence $A \vee(S^{n})^{\vee p} \to B \vee (S^{n})^{\vee q}$, for some $p,q \in \N$.  In fact this simple equivalence can be chosen so that it agrees with $f$ over $A^{n-1}$. This is  Theorem 14 of \cite{W5}.
\end{Remark}

\begin{Proof} {\bf (Theorem \ref{claim1})} Consider the filtered space:
  $$P=(P,P_L,P_{L-1},...,P_0=*),$$
 where  $P$ is the reduced mapping cylinder of $F\colon M \to N$ (which has a  natural CW-decomposition), and $P_k$ is the reduced mapping cylinder of $F^k=F_{|M^k}\colon M^k \to N^k$, for $k=1,...,L-1$. We  further define:
 \[P_L=M^L \cup P_{L-1} \cup N^L=P^L.\]
Therefore each $P_k$ is embedded cellularly in $P_{k+1}$ and $P$, for $k=0,...,L-1$.

We want to prove that there exist filtered homotopy equivalences:
\[ \bar{M} \cong P \cong \bar{N}. \]

By Lemma  \ref{retraction}, there exists a deformation retraction of  $P$ onto $M$, say $\r\colon P \times I \to P$,  which we can suppose to be cellular, such   that the restriction of $\r$ to $P_k\times I$ realises a  deformation retraction of $P_k$ in $M^k$, for $k=0,...,L-1$. Define: \[r(u)=\r(u,1), u \in P.\] Therefore $r$ is a cellular  map $P \to M$.

Let $\{e^L_1,...,e^L_{n_L}\}$ be the set of $L$-cells of $N$. Each of these
cells attaches to $N^{L-1}$ along a certain map: 
\[\f^L_{i_L}\colon  S^{L-1}_{i_L} \doteq \d(e^L_{i_L})\to N^{L-1},i_L=1,....,n_L.\]  One obtains $P_L$ from $M^L \cup P_{L-1}$ by attaching these $L$-cells to it  in the obvious way. Notice that $N^{L-1}$ is included in $P_{L-1}$, cellularly.

The deformation retraction $\r$, from $P$ onto $M$,  gives us a homotopy (in
$P_{L-1}$) connecting $\f\il\colon S^{L-1}_{i_L} \to N^{L-1}\subset P_{L-1}$
with $(r\circ \f\il)\colon S^{L-1}_{i_L}\to M^{L-1}\subset P_{L-1}$, for
$i_L=1,....,n_L$, since the restriction of $\r$ to $P_{L-1} \times I$ is a
deformation retraction from $P_{L-1}$ onto $M^{L-1}$. On the other hand, the
restriction $r_{|e^L_{i_L}}$ of $r$ to the $L$-cell $e\il$ will provide a null
homotopy $(r \circ \f\il) \to *$ (in $M^L$) for $i_L=1,....,n_L$ (recall that
$r$ is cellular, thus $r  (e^L_{i_L})\subset M^L$). Note that we denote by $*$
both the base point and any function with values in $\{*\}$, and the same for
any singleton $\{a\}$. We have proved that  each  attaching map $\f\il\colon  \d(e \il) \to M^L \cup P_{L-1}$, where $i_L=1,...,n_L$, is null homotopic.

Note that the CW-complex  $M^L$ is a deformation retract of $ M^L \cup P_{L-1}$, since                
$M^{L-1} \subset M^L$ is a deformation retract of $ P_{L-1}$, and
$M^{L-1}=P_{L-1} \cap M^L$. In particular the inclusion map $M^L \ra{i} M^L
\cup P_{L-1}$ is a homotopy equivalence. Therefore, there exist homotopy
equivalences:
\[(M^L,*) \vee \V_{i_L=1}^{n_L} (S^L_{i_L},*) \ra[\cong]{i} (M^L \cup
P_{L-1},*)\vee \V_{i_L=1}^{n_L} (S^L_{i_L},*) \ra[\cong]{} (P_L,*),\] 
the second one implied by the fact that the attaching maps
$\f\il,i_L=1,...,n_L$ are null homotopic. 
Let $G_L$ be their composition, thus  $G_L$ defines a homotopy equivalence: 
\[G_L\colon (M^L,*)\vee \V_{i_L=1}^{n_L} (S^L_{i_L},*) \to  (P_L,*).\]

Note that the restriction of $G_L$ to $M^L$ is (or can be taken to be) the inclusion map $M^L \to P_L$.  By Corollary \ref{retraction}, the restriction of  this inclusion map to $M^{L-1}$ defines a filtered homotopy equivalence:
$(M^{L-1},....,M^0=*) \cong (P_{L-1},...,P_0=*)$. 
 Therefore it follows, from Theorem \ref{Simplify}, that $G_L$ is also  a filtered homotopy equivalence:
\begin{multline*} 
G_L\colon (M^L,M^{L-1},....,M^0=*) \vee \V_{i_L=1}^{n_L} (S^L_{i_L},*,..,*)\\\to (P_L,P_{L-1},...,P_0=*). 
\end{multline*} 

Given that  the restriction of $G_L$ to $M^L$ is  the restriction of the inclusion map $M \to P$ to $M^L$, we can conclude that the map $G_L$ extends to a filtered map (which we also call $G_L$):
\begin{multline*} 
G_L\colon (M,M^L,M^{L-1},....,M^0=*) \vee \V_{i_L=1}^{n_L} (S^L_{i_L},S^L_{i_L},*,..,*)\\\to (P,P_L,P_{L-1},...,P_0=*),
\end{multline*} 
which  when restricted to $M$  is the inclusion map  $M \to P$.

 We want to prove that $G_L$ extends to a map $G: M \vee \left (D^{L+1}
 \right)^{\vee n_L}= \bar{M} \to P$. To this end, we need to prove that the restriction of $G_L$ to each $S^L_{i_L} \cong D^L_{i_L}/\d D^L_{i_L}$ is null homotopic.

Let us be more explicit about the construction of $G_L$. It is important to remember that  the map $G_L$ depends, \emph{explicitly}, on the null  homotopies: 
 \[H_{i_L}\colon * \to \left ( \f\il\colon \d(e\il) \to N^{L-1} \subset M^{L} \cup P_{L-1}\right ),\]
 which we choose. Here $i_L=1,...,n_L$. These homotopies are taken in  $M^L \cup P_{L-1}$.  Let us describe each $H_{i_L}\colon \d(e\il)\times I  \to M^L \cup P_{L-1} $  in full detail: 

Let $i_L\in \{1,...,n_L\}$. Recall $S^{L-1}_{i_L}\doteq \d(e\il)$. Define $H_{i_L}\colon  S^{L-1}_{i_L} \times I \to M^L \cup P_{L-1} $ in the following way: Let $a^L_{i_L}$ be the central point   of the closure  $\overline{e^L_{i_L}}\cong D^L$ of $e\il$. Therefore, there exists a path $\g^L_{i_L}$ in $M^L$ connecting  the base point $*$ with  $r(a^L_{i_L})$. When $t \in [0,1/3] \subset I$, we define $H_{i_L}$ directly from $\g^L_{i_L}$, in the obvious way. This yields a homotopy $J^1_{i_L}\colon  * \to r(a^L_{i_L})$, in $M^L$.
Specifically:

\[J^1_{i_L}(x,t)=\g \il(t), \forall x \in S^{L-1}_{i_L}, \forall t \in I.\]

When $t \in [1/3,2/3] \subset I$, we consider $H_{i_L}$ to be defined, in the
natural way, from the restriction $r_{|e^L_{i_L}}$ of $r$ to $e\il$. This
defines a homotopy $J^2_{i_L}\colon r(a^L_{i_L}) \to (r \circ \f^L_{i_L})$, in
$M^L$. Explicitly (note $\overline{e\il} \cong D^L=[-1,1]^L$):

\[J^2_{i_L}(x,t)=r_{|e\il}(xt), \forall x \in S^{L-1}_{i_L}, \forall t \in  I=[0,1].\]
Recall that $e\il$ attaches along $\f\il\colon S^{L-1}_{i_L} \to N^{L-1}\subset M^L \cup P_{L-1}$.

 Finally, the deformation retraction $\r\colon P \times I \to P$ (of $P$ in
$M$) will define a homotopy $J^3_{i_L}$ (in $P_{L-1}$) connecting  $(r \circ
\f^L_{i_L})\colon S^{L-1}_{i_L} \to M^{L-1} \subset M^L \cup P_{L-1}$ with $\f^L_{i_L}\colon S^{L-1}_{i_L} \to N^{L-1} \subset M^L \cup P_{L-1}$. Explicitly:

\[J^3_{i_L}(x,t)=\r(\f\il(x),1-t), \forall x \in S^{L-1}_{i_L}, \forall t \in I.\]
 It is important to remember that the restriction of $\r$ to $P_{L-1} \times
 I$ defines a deformation retraction of $P_{L-1}$ in  $M^{L-1}$. This is an
 extremely crucial fact. 

 We thus define $H_{i_L}\colon * \to \f^L_{i_L}$ as the concatenation: \[H_{i_L}=\left (J^1_{i_L}J^2_{i_L}J^3_{i_L}\right)\colon S^{L-1}_{i_L} \times I \to M^L \cup P_{L-1},\]
 for $i_L=1,...,n_L$.

  As we have seen, when restricted to $M^L$, the map $G_L$ is the restriction of the inclusion map $M \to P$ to $M^L$. 
For each $i_L=1,...,n_L$, the restriction $K\il$ of $G_L$ to $S^L_{i_L}=D^L_{i_L} /\d D^L_{i_L}$ is given, in the obvious way, from\footnote{This is the standard homotopy equivalence $A \cup_f D^n \to B \cup_g D^n$, where $A$ is some space, constructed from a homotopy   $H\colon S^{n-1} \times I \to A$ connecting $f\colon S^{n-1} \to A$ with $g\colon S^{n-1} \to A$. See for example \cite{H}, proof of  Proposition 0.18.}: 

\begin{enumerate}
\item The homotopy, $H_{i_L}\colon S^{L-1}_{i_L} \times I \to M^L \cup
  P_{L-1}$ (connecting $*$ with the map $\f^L_{i_L}\colon S^{L-1}_{i_L} \to
  N^{L-1}\subset  M^L \cup P_{L-1}$), for $1/3\leq |x| \leq 1, x \in D^L_{i_L} /\d D^L_{i_L} $.
\item The obvious re-scaling of the restriction  $\id_{|e\il}$ of the identity map $\id\colon N \to N$ to the cell $e^L_{i_L}
$, for   $0\leq |x| \leq 1/3, x \in D^L_{i_L} /\d D^L_{i_L}$. (Recall that the
cell $e\il$ attaches along $\f^L_{i_L}\colon S^{L-1}_{i_L} \to N^{L-1}$).
\end{enumerate} 
 We will make precise the construction of $K\il$ is a second. The main point is that when $K\il={G_L}_{|S^L_{i_L}}$ is  considered to be  a map $S^L_{i_L}=D^L_{i_L} /\d D^L_{i_L}\to P$, it is homotopic (modulo the border) to the map $K'\il$ which is defined, in the natural way, from: 
\begin{enumerate}
\item the homotopy $J^1_L J^2_L$ (connecting $*$ with $(r \circ \f^L_{i_L})$), for $2/3\leq |x| \leq 1$ , 
\item the constant homotopy $ (r \circ \f^L_{i_L}) \to (r \circ \f^L_{i_L}) $, for $1/3 \leq |x| \leq 2/3$,
\item the obvious re-scaling of the function $r_{|e^L_{i_L}}$, for  $0\leq |x| \leq 1/3$. (Recall again that $e\il$ attaches along $\f^L_{i_L}\colon S^{L-1}_{i_L} \to N^{L-1}$).
\end{enumerate}
Indeed, a homotopy $W\il\colon D^L_{i_L}/D\il\times I \to M^L \cup P_{L-1}$ connecting $K\il$ with $K'\il$ is given by:

\begin{multline*}
W\il(x,t)=\\\left \{ \begin{CD} J^1_{i_L}J^2_{i_L}(x/|x|,3(1-|x|)), \textrm{ if } 2/3 \leq |x| \leq 1, \forall t \in I,  \\ 
             \r_1 \left (\f\il(x/|x|),t+ (3-3t) (|x|-1/3)\right)\textrm{ if } 1/3 \leq |x| \leq 2/3, \forall t \in I , 
   \\\r_1\left (\id_{|e\il}(3x),t\right), \textrm{ if } 0 \leq |x| \leq 1/3, \forall t \in I.  \\ 
          \end{CD} \right .
\end{multline*} 
This also gives an explicit  definition of $K\il(x)=W\il(x,0)$ and $K'\il(x)=W\il(x,1)$, where $x \in D^L_{i_L}/\d D^L_{i_L}$.

The map $K'\il\colon  D^L_{i_L}/D^L_{i_L} \to M^L \subset M^L \cup P_{L-1}$ is, obviously, null homotopic. We have, therefore, proved that $G_L$ extends to a  map: 
\begin{multline*} 
G\colon \bar{M}\doteq (M,M^L,M^{L-1},....,M^0=*) \vee \V_{i_L=1}^{n_L} (D^{L+1}_{i_L},S^L_{i_L},*,..,*)\\\to (P,P_L,P_{L-1},...,P_0=*). 
\end{multline*}

Let us prove that $G$ is a (filtered) homotopy equivalence. Let $G^k$ be the restriction of $G$ to $\bar{M}^k$, for $k=0,...,L$. As we have seen, for $k=0,...,L-1$, the map $G^k\colon M^k=\bar{M}^k \to P_k$ is a homotopy equivalence. By construction, $G^L$ is a homotopy equivalence as well. Obviously $G$ is also a homotopy equivalence, since its restriction to $M$ is a homotopy equivalence. Therefore, by Theorem \ref{Simplify}, it follows that $G$ is a filtered homotopy equivalence.

Repeating the same argument for the inclusion of $N$ in $P$, we can prove that
there exists a filtered homotopy equivalence:
\begin{multline*} 
Q\colon \bar{N}\doteq(N,N^L,N^{L-1},....,N^0=*) \vee \V_{j_L=1}^{m_L} (D^{L+1}_{j_L},S^L_{j_L},*,..,*)\\\to (P,P_L,P_{L-1},...,P_0=*). 
\end{multline*}
Here $m_L$ is the number of cells of $M$ of order $L$. 
Note that both $G$ and $Q$, when restricted to $N^{L-1}$ and $M^{L-1}$ are the inclusion maps in $P$. 

The obvious retraction $R\colon P_{L-1} \to N^{L-1}=\bar{N}^{L-1}$ (see the discussion after Corollary \ref{retraction}) is a filtered homotopy inverse of the map $Q^{L-1}\colon  (N^{L-1},...,N^0)\\\to (P_{L-1},...,P_0)$. Therefore, there exists a filtered homotopy inverse $R'$ of $Q$ extending $R$, by Corollary \ref{strenght}. Note that $(R \circ G)_{|M^{L-1}}=F_{|M^{L-1}}.$
In particular there exists a filtered homotopy equivalence:
\begin{multline*} 
F'=R' \circ G\colon \bar{M}=(M,M^L,M^{L-1},....,M^0=*) \vee \V_{i_L=1}^{n_L} (D^{L+1}_{i_L},S^L_{i_L},*,..,*)\\
\to (N,N^L,N^{L-1},....,N^0=*) \vee \V_{j_L=1}^{m_L} (D^{L+1}_{j_L},S^L_{j_L},*,..,*)=\bar{N},
\end{multline*}
with $F'_{|M^{L-1}}=F_{|M^{L-1}}.$
We now need to prove that $F'$ can be chosen to be  cellular. By definition, $F'_ {|\bar{M}^L }\colon  \bar{M}^L \to \bar{N}^L$ is cellular. The map $F'$ is homotopic (modulo   $\bar{M}^L$) to a cellular map  $F''$, thus $F''$ is a homotopy equivalence.  Since the restriction of  $F''$ to each $\bar{M}^k$ is a homotopy equivalence for $k=0,...,L$, it follows that $F''\colon  \bar{M} \to \bar{N}$ is a filtered homotopy equivalence, by Theorem \ref{Simplify}. The map $F''$ agrees with $F$ over $M^{L-1} =\bar{M}^{L-1}$.  This finishes the proof of Theorem \ref{claim1}.
\end{Proof}

\section{Applications to the Fundamental Crossed Complex of a CW-Complex}

\subsection{The Dependence of  $\Pi(M)$ on the Cell Decomposition of $M$}
Let $M$ be a CW-complex, equipped with its skeletal filtration. As we have
mentioned before in \ref{Homotopy}, even though the crossed complex $\Pi(M)$ depends strongly on the cellular decomposition of $M$, J.H.C. Whitehead proved that
the homotopy type of $\Pi(M)$ depends only on the homotopy type of $M$, see
\cite{W4}. See also Remark \ref{simplerefer}, below.  The following result  (Theorem \ref{Main2}) strengthens this 
slightly.  We freely use the results and notation of \ref{ccomplex}.

If $n \in \N$, define $\D^n=(\D^n_m,\d_m)$ as the crossed complex, 
which is the trivial group for dimensions not equal to  $n$ or $n-1$, whilst
$\D^n_n=\Z$ and $\D^n_{n-1}=\Z$, the border map $\d_n\colon \Z \to\Z$ being the
identity. For $n=2$,  the action of $\Z$ in $\Z$ is defined to be the trivial action. 
 Therefore $\D^n\cong \Pi(D^n)$, where the $n$-disk $D^n$ is given its natural
CW-structure with one 0-cell, one $(n-1)$-cell, and one $n$-cell.

 Let  $M$ be a $CW$-complex, with a unique 0-cell. Recall that, by definition,  $\Pi_L(M)=\Pi_L(M,M^{L-1},M^{L-2},...,M^0=*)$.

\begin{Theorem}\label{Main2}
Let $M$ and $N$ be  CW-complexes, both with a unique 0-cell which we take to
be their base points. Give $M$ and $N$  their skeletal filtrations.  Suppose
that $M$ and $N$ are homotopic as topological spaces. There exists an isomorphism of crossed complexes:
\begin{equation*}
\Pi(M) \V_{n=1}^{\infty} \left ((\D^{n+1})^{\vee \left (l^M_{\{M,N\}}(n)\right) }\right)\cong \Pi(N) \V_{n=1}^{\infty}\left ( (\D^{n+1})^{\vee \left (l^N_{\{M,N\}}(n)\right) }\right).
\end{equation*}
Moreover, if $L\in \N$ we also have 
\begin{equation*}
\Pi_L(M) \V_{n=1}^{L-1} \left ((\D^{n+1})^{\vee \left (l^M_{\{M,N\}}(n)\right) }\right)\cong \Pi_L(N) \V_{n=1}^{L-1}\left ( (\D^{n+1})^{\vee \left (l^N_{\{M,N\}}(n)\right) }\right).
\end{equation*}
\end{Theorem}
The constants $l^N_{\{M,N\}}(n)$ and $  l^M_{\{M,N\}}(n)$, where $n \in \N$,
were defined at the start  of \ref{general}. They also make sense, though may
be  infinite,  if $M$ or $N$ have infinite cells of order $L$, for some $L \in \N.$

\begin{Proof}
This follows immediately from Theorem \ref{Main} and (R. Brown and
P.J. Higgins') Corollary \ref{Free}. 
\end{Proof}

One simple consequence  of this theorem is a new proof of (J.H.C. Whitehead's) Theorem \ref{homotopy}:

\begin{Corollary}
Let $M$ and $N$ be  CW-complexes with unique 0-cells. Suppose that they are homotopic as (based) spaces. There exists a homotopy equivalence $\Pi(M) \cong \Pi(N)$. In fact, for any $L \in \N$, the crossed complexes $\Pi_L(M)$ and $\Pi_L(N)$ are also homotopy equivalent. 
\end{Corollary}
\begin{Proof}
This result is an easy consequence of Corollary \ref{extension} and the
previous theorem. Also needed  is the universal property defining free
products of crossed complexes, as well as Lemma \ref{connection}.
\end{Proof}

\begin{Remark}\label{simplerefer}
Theorem \ref{Main2} has a large intersection with Theorem 16 of \cite{W5}. Roughly
speaking, this result due to J.H.C. Whitehead  says that if $f\colon \A \to \B$ is a
homomorphism of finite, $L$-truncated, totally free crossed complexes  inducing an isomorphism on the
homotopy and homology groups of $\A$ and $\B$ up to dimension $L-1$, then it follows that there
exists a simple homotopy equivalence $f'\colon \A'\supset \A \to \B'\supset \B$, which when  restricted
to $\A\subset \A'$ is  $f$. Here $\A'$ and $\B'$
are obtained from $\A$ and $\B$ by making a certain number of free products with $\Pi_L(S^L)$,
where the $L$-sphere $S^L$ is given a cell decomposition with unique cells of order $0$ and
$L$. The notion of simple homotopy equivalence of totally free crossed
complexes is also defined in \cite{W5}. It is very similar with 
the relation between crossed complexes indicated in Theorem \ref{Main2}. 
\end{Remark}
\subsection{A Counting Type Homotopy Invariant $I_\A$}
\subsubsection{The Definition of $I_\A$}
We now need to restrict our discussion to CW-complexes which only have a finite
number of $L$-cells for each $L \in \N$, therefore avoiding infinities. As
usual, all CW-complexes that  we consider have  unique 0-cells, taken to be their base points.

\begin{Definition}
A (reduced) crossed complex $\A=(A_n)$ is called finite if $\A$ is $L$-truncated for some $L$ and, moreover, all groups $A_n, n \in \N$ are finite. 
\end{Definition}

Let $\A=(A_n)$ be a finite $L$-truncated crossed complex. Therefore, if $M$ is
a CW-complex with a finite number of cells of each dimension $L \in \N$, then the
number of morphisms $\Pi(M) \to \A$ is finite. This follows from Lemma \ref{colourings}. Moreover, by Proposition \ref{Cotrunc}, there exists a one-to-one correspondence between morphisms $\Pi(M) \to \A$ and morphisms $\Pi_L(M) \to \A$.

For a set $A$, denote its cardinality by $\#(A)$. As usual if $M$ is a CW-complex, we denote the number of cells of $M$ of dimension  $n$ by $l^M_n$.

\begin{Theorem}\label{counting}
Let $M$ be a crossed complex with a unique $0$-cell, such that $M$ only has a
finite number of cells in each dimension (or, alternatively, up to dimension
$L$). Let also $\A=(A_n)$ be an $L$-truncated finite crossed complex. Define:
\begin{align*}
I_\A(M)&=\# \left(\Hom(\Pi(M),\A)\right) \prod_{n=1}^{\infty} \left (
\prod_{m=1}^\infty  \#\left (A_{m+n}\right)^{l^M_m} \right )^{(-1)^n}\\
      &=\# \left(\Hom(\Pi_L(M),\A)\right) \prod_{n=1}^{\infty} \left (
\prod_{m=1}^\infty  \#\left (A_{m+n}\right)^{l^M_m} \right )^{(-1)^n}.
\end{align*}
Then $I_\A(M)$ does not depend on the CW-decomposition of $M$, and it is a homotopy invariant of $M$.
\end{Theorem}
\begin{Remark}Note that Lemma \ref{colourings} ensures that $I_\A(M)$ is, in
  principle, calculable, in a combinatorial way.
\end{Remark}

Theorem \ref{counting} is a consequence of the following lemma, easy to prove:
\begin{Lemma}
Let $\A=(A_n)$ be a crossed complex. There exists a one-to-one correspondence between morphisms $\D^n \to \A$ and elements of $A_n$, where $n \in \N$. 
\end{Lemma}
Recall $\D^n=\Pi(D^n)$, where $D^n$ is provided with its natural cell decomposition with unique cells of order $0$, $(n-1)$ and $n$.

\begin{Proof} {\bf (Theorem \ref{counting}) }
Let $M$ and $N$ be homotopic cellular spaces. Give $M$ and $N$
CW-decompositions with a unique $0$-cell, and such that $M$ and $N$ only have
a finite number of $L$-cells for each $L \in \N$. By Theorem \ref{Main2} we
can conclude:
\begin{multline*}
\# \Hom\left (\Pi(M) \V_{n=1}^{\infty} \left ((\D^{n+1})^{\vee \left (l^M_{\{M,N\}}(n)\right) }\right),\A \right) \\= \#\Hom  \left( \Pi(N) \V_{n=1}^{\infty}\left ( (\D^{n+1})^{\vee \left (l^N_{\{M,N\}}(n)\right) }\right), \A \right).
\end{multline*}
Therefore, from the universal property of free products of crossed complexes, as well as the previous lemma, it follows that:
\begin{multline*}
\#\left ( \Hom (\Pi(M), \A)\right) \prod_{n=1}^{\infty} \#(A_{n+1})^{l^M_{\{M,N\}}(n)}\\=\# \left (\Hom (\Pi(N), \A)\right) \prod_{n=1}^{\infty} \#(A_{n+1})^{l^N_{\{M,N\}}(n)}.
\end{multline*}
The result follows from some straightforward algebra.
\end{Proof}

\subsubsection{Geometric Interpretation of $I_\A$}
The results  described in  \ref {Homotopy}, \ref{conthom} and \ref{Classifying}, as well as the notation introduced,  will be used actively.

Let $M$ be a path connected space. Suppose  that $M$ only has a finite number of
non-trivial homotopy groups, all of which are finite. Cellular spaces of this type are
studied in  \cite{E,L}, for example. They generalise Eilenberg-McLane spaces
of finite groups.  Define the following ``multiplicative Euler Characteristic type'' invariant:

\[\X(M)= \prod_{k=1}^{\infty} \left [ \# (\pi_k(M)) \right ]^{(-1)^k}.\]
In general, if $M$ has a finite number of connected components, define:
\[\X(M)= \sum_{M_0 \in \pi_0(M)} \quad  \prod_{k=1}^{\infty} \left [\#
  (\pi_k(M_0))\right ]^{(-1)^k},\]
where the sum is extended to all connected components $M_0$ of $M$.

Let $\A$ be a finite (non-reduced) crossed complex, as usual in the sense that
$\A$ is $L$-truncated, for some $L \in \N$, and all the groupoids $A_n,n \in \N$ are finite. Then
the classifying space $|\A|$ of $\A$ is a topological space with a finite
number of homotopy groups, all of which are finite, by  (R. Brown and P.J. Higgins') Theorem \ref{homgroups}.  We warn the reader that  classifying spaces of crossed complexes do not exhaust all spaces of this type. Let $\B$ be a (reduced)
finite 
crossed complex. If $M$ is a CW-complex with a finite number of $L$-cells for any $L \in \N$, and a unique 0-cell, then the (non-reduced) crossed complex
$CRS(\Pi(M),\B)$ is finite, by Lemma \ref{colourings} and Corollary \ref{cardinality}. In particular,  the  quantity:
\[\X(|CRS(\Pi(M),\B)|)=\X\left (TOP\left ((M,*),(|\B|,*)\right ) \right ),\]
is finite. Recall (R. Brown and P.J. Higgins') Theorem \ref{weaktype}.

\begin{Lemma}
Let $\A=(A_n)$ be a (non reduced) finite crossed complex. Let $C$ be the object set
of $A_1$. As usual, we denote by  $A_n^c$, the set of morphisms of $A_n$ with source $c$, where   $n \in N$ and $c \in C$. We have:
\[\X(|\A|)=\sum_{c \in C} \prod_{k=1}^{\infty} [\#(A^c_k)]^{(-1)^k}.\] 
\end{Lemma}
\begin{Proof}
Let $\pi_0(A_1)$ be the set of connected components of the groupoid $A_1$. There exists a
 one-to-one correspondence between $\pi_0(A_1)$ and the set $\pi_0(|\A|)$ of connected
 components  of $|\A|$. Each element of $c \in C$ yields a unique 0-cell $c$
 of $|\A|$. Let also $[c]$ denote the connected component of $A_1$ to which 
$c$ belongs.  Recall that $\d^c_n$ equals $\d_n$ restricted to ${A_n^c}$, if $c
 \in C$ and $n \in \N$. We have:
\begin{align*}
\X(|\A|)&= \sum_{M \in \pi_0(|\A|)} \prod_{k=1}^{\infty} \left ( \# (\pi_k(M))
         \right )^{(-1)^k}\\
         &=\sum_{c \in C} \frac{1}{\#([c])} \prod_{k=1}^{\infty}
         \left ( \# (\pi_k(|\A|,c)) \right )^{(-1)^k}\\
         &=\sum_{c \in C} \frac{1}{\#(\pi_1(|\A|,c)) \#([c])} \prod_{k=2}^{\infty}
         \left (\frac{\#(\ker (\d^c_k))} {\#({\rm im}(\d^c_{k+1}))}
         \right)^{(-1)^k}\\
         &=\sum_{c \in C} \frac{1}{\#(\pi_1(|\A|,c)) \#([c])} \prod_{k=2}^{\infty}
         \left (\frac{\#(\ker (\d^c_k)) \#(\ker (\d^c_{k+1}))} {\#(A_{k+1}^c)}
         \right)^{(-1)^k}\\
          &=\sum_{c \in C} \frac{\#(\ker (\d^c_2))}{\#(\pi_1(|\A|,c)) \#([c])} \prod_{k=2}^{\infty}
         \left ({\#(A_{k+1}^c)} \right)^{(-1)^{(k+1)}}.\\
\end{align*}
If $c,d \in C$, recall that  $A^{(c,d)}_1$ denotes the set of morphisms of $A_1$ with  source $c$ and target $d$. We
have $\pi_1(|\A|,c)=A_1^{(c,c)}/{\rm im} ( \d^c_2)$. On the other hand, since $A_1$ is a groupoid, it follows that
$\#(A_1^{(c,c)})\#([c])=\#(A^c_1), \forall c \in C$. In particular we have:
\begin{align*}
\X(|\A|)&=\sum_{c \in C} \frac{\#(\ker (\d^c_2)) \#({\rm im} ( \d^c_2))} 
{\#(A^c_1)} \prod_{k=2}^{\infty}
         \left ({\#(A_{k+1}^c)} \right)^{(-1)^{(k+1)}}\\
         &=\sum_{c \in C} \frac{\#(A^c_2)} 
{\#(A^c_1)} \prod_{k=2}^{\infty}
         \left ({\#(A_{k+1}^c)} \right)^{(-1)^{(k+1)}}\\
         & =\sum_{c \in C} \prod_{k=1}^{\infty}
         \left ({\#(A_{k}^c)} \right)^{(-1)^{k}}.\\
\end{align*}
\end{Proof} 

From Corollary \ref{cardinality} we can deduce the following theorem, which  at the same time  interprets and gives an alternative proof of the existence of the invariant $I_\A$, where $\A$ is a finite crossed complex.

\begin{Theorem}
Let $\A$ be a finite crossed complex. Let $M$ be a CW-complex with a unique 0-cell and a finite number of $L$-cells for any $L \in \N$. We have:
\[I_\A(M)=\X \left (TOP((M,*),(|\A|,*))\right).\]
\end{Theorem}
\begin{Proof}
As we have mentioned before (Theorem \ref{weaktype}), there exists a map: 
\[\psi\colon |CRS(\Pi(M),\A)| \to TOP((M,*),(|\A|,*)) ,\]
which is  a weak homotopy equivalence; a  result due to R. Brown and P.J. Higgins, appearing in \cite{BH5}. Therefore:
\[
\X \left (TOP((M,*),(|\A|,*))\right)= \X \left(|CRS(\Pi(M),\A)|\right).\]
We now need to apply the previous lemma together with Corollary \ref{cardinality}.
\end{Proof}
\begin{Remark}
These final results originated from  discussions that  I had with Ronnie Brown and Tim
Porter. In fact the argument for the simpler case of crossed modules was
initiated by them.
\end{Remark}

\subsubsection{ A Short Discussion and an Extension of $I_\A$}
Let $\A=(A_n,\d_n)$ be a finite crossed complex.
In \cite{Y}, D.  Yetter defined a 3-manifold invariant for any finite crossed
module. The invariant $I_\A$ explains Yetter's invariant, in the cellular category,  and
generalises it to crossed complexes.  A previous systematic study and
upgrading of Yetter's construction was done by T. Porter, and  appeared in
\cite{P1,P2}. On the other hand, the  construction due
to M. Mackaay of
4-manifold invariants  which appears in \cite{Mk} is, conjecturally, related
with 3-types, having a further incorporation of cohomology classes of them in
the manner shown below.

The crossed module case was also studied  in \cite{FM1,FM2}. The first article
considered only  knot complements, whereas  the second article covered general
CW-complexes (for example subsection \ref{two} is almost extracted from \cite{FM2}). One
of the main conclusions of this work was that, in the case of knotted
surfaces, the invariant $I_\G$ (where $\G$ is a finite crossed module) is a non-trivial, very calculable invariant.  In fact \cite{FM2} contained an algorithm for the  calculation of $I_\G$ from  movie presentations of knotted embedded surfaces in $S^4$.

Let $M$ be a finite CW-complex. We can easily describe a morphism $\Pi(M) \to
\A$.  Recall  Lemma \ref{colourings}. Roughly speaking,  these morphisms are
specified, uniquely,  by their value on the $L$-cells of $M$, as long as we
choose for any $L\in \N$ and any   $L$-cell $e^L$ 
an element of $\pi_{L-1}(M^{L-1},*)$ (defined up to acting by an element of
$\pi_1(M^1,*)$) along which $e^L$ attaches to $M^{L-1}$. After making these
choices, there exists a one-to-one correspondence between  morphisms $\Pi(M) \to
\A$ and   colourings of each $L$-cell of $M$ ($L \in \N$) by an element of $A_L$, with the
obvious compatibility relations with the boundary maps of $\Pi(M)$ and $\A$. Therefore,  it is  not a difficult task to calculate $I_\A(M)$ where $M$ is a CW-complex with a
unique 0-cell, and its calculation is of a  combinatorial nature, similar to
the calculation of the cellular homology groups of $M$. We refer to \cite{FM2} for some calculations in the crossed module case.

The invariant $I_\A$ can be  naturally  twisted by $n$-dimensional cohomology
classes $\w$ of $|\A|$ (the classifying space of $\A$) to an invariant $I_M(\quad,\w)$, as long as we restrict $I_\A(\quad,\w)$ to oriented $n$-dimensional closed
manifolds. Indeed, if $\w \in H^n(|A|)$, we can define: 
\begin{align*}
I_\A(M,\w)&=\sum_{f \in [(M,*), (|\A|,*)]} \left <o_M,f^* (\w) \right>\\ 
&\quad \quad \quad \quad\quad\quad\quad\quad\prod_{k=1}^{\infty}  \#(\pi_k(TOP((M,*),(|\A|,*)),f))^{{(-1)}^k}\\
          &= \sum_{F \in \Hom (\Pi(M),\A)   } \left <o_M,F_g^* (\w) \right> \prod_{n=1}^{\infty} \left (
\prod_{m=1}^\infty  \#\left (A_{m+n}\right)^{l^M_m} \right )^{(-1)^n},
\end{align*}
where $o_M$ is the orientation class of $M$, and, as usual, $l^M_n$ denotes the
number of cells of $M$ of order $n \in \N$.  In addition, if
$F \in \Hom(\Pi(M),\A)$, then $F_g:(M,*) \to (|A|,*)$ denotes a geometric
realisation of $F$, which is uniquely defined up to homotopy. See \cite{BH5}. The algebraic description of the cohomology of the classifying space of crossed modules (particular cases of crossed complexes) appears for example in \cite{E2,Pa}.

From its very definition, the invariant $I_\A(\quad,\w)$ is a generalisation of the Dijkgraaf-Witten invariant of manifolds (see \cite{DW}). If $M$ is provided a triangulation, then $I_\A(\quad,\w)$ can   be calculated in a 
similar form to the Dijkgraaf-Witten,  but considering colourings on
any simplex of $M$, rather than only on the edges of $M$, at least in the untwisted case. We will consider these issues in a subsequent
publication. 

\section*{Acknowledgements}

This work was financed by  Funda\c{c}\~{a}o para a Ci\^{e}ncia e Tecnologia (Portugal), post-doctoral grant number SFRH/BPD/17552/2004, part of the research project POCTI/MAT/60352/2004 (``Quantum Topology''), also financed by F.C.T.

I would like to express my gratitude to Gustavo Granja for many helpful discussions, which were extremely important for the  development of this work.  
  I also  want to thank Tim Porter and Ronnie Brown for introducing me to this
  subject, and  for a large amount of  further suggestions and help, which influenced this article, and finally Roger Picken for his constant support.

\end{document}